\input ./style/arxiv-general.cfg
\documentclass[aos,MSNbibl,nameyear,seceqn,dvips]{arximspdf}
\makeatletter
   \@ifpackageloaded{graphicx}{}{\usepackage{graphicx}}
\makeatother
\doi{10.1214/15-AOS1323}% Updated by VTEXPTS2LaTeX.exe, 27.08.2015 13:08
\volume{44}
\issue{1}
\pubyear{2016}
\firstpage{1}
\lastpage{30}
\docsubty{FLA}

\makeatletter
\def\cal{\mathcal}
\newcommand{\eqref}[1]{(\ref{#1})}
\newtheorem{thmm}{Theorem}

\newtheorem{prop}{Proposition}

\newproclaim{assumption}{Assumption}
\newproclaim{remark}{Remark}
\newproclaim{definition}{Definition}
\makeatother

\begin{document}
\begin{frontmatter}

\title{Functional linear regression with points of impact}%
%\thanksref{T1}}
%\thankstext{T1}{Financial supported by the DFG through GRK 1707 and
%the Hausdorff Center for Mathematics is greatly acknowledged}
\runtitle{Point impact}

\begin{aug}
% Corresponding author: Alois Kneip - akneip@uni-bonn.de% Updated by
%VTEXPTS2LaTeX.exe, 27.08.2015 13:08
\author[A]{\fnms{Alois}~\snm{Kneip}\corref{}\thanksref{T1}\ead[label=e1]{akneip@uni-bonn.de}},
\author[B]{\fnms{Dominik}~\snm{Po{\ss}}\thanksref{T2}\ead[label=e3]{dposs@uni-bonn.de}}
\and
\author[C]{\fnms{Pascal}~\snm{Sarda}\ead[label=e2]{Pascal.Sarda@math.ups-tlse.fr}}
\runauthor{A. Kneip, D. Poss and P. Sarda}

\thankstext{T1}{Supported by the DFG through GRK 1707 and the Hausdorff
Center for Mathematics.}% is greatly acknowledged}
\thankstext{T2}{Supported by the DFG through GRK 1707.} % and the
%Hausdorff Center for Mathematics is greatly acknowledged}
\affiliation{Universit\"{a}t Bonn and Institut de Math\'ematiques de Toulouse}
\address[A]{A. Kneip\\
Institut f\"ur Finanzmarkt\"okonomik und Statistik\\
Department of Economics and Hausdorff\\
Center for Mathematics\\
Universit\"at Bonn\\
Adenauerallee 24-26\\
53113 Bonn\\
Germany\\
\printead{e1}}

\address[B]{D. Poss\\
Bonn Graduate School of Economics\\
Department of Economics\\
Institut f\"ur Finanzmarkt\"okonomik\\
\quad und Statistik \\
Universit\"at Bonn\\
Adenauerallee 24-26\\
53113 Bonn\\
Germany\\
\printead{e3}}

\address[C]{P. Sarda\\
Institut de Math\'ematiques de Toulouse; UMR 5219 \\
Universit\'{e} de Toulouse and CNRS\\
118, Route de Narbonne\\
31062 Toulouse Cedex\\
France\\
\printead{e2}}
%\runauthor{}
%\dedicated{}
\end{aug}

% HISTORY:
%
\received{\smonth{7} \syear{2014}}% Updated by VTEXPTS2LaTeX.exe,
%27.08.2015 13:08
%
\revised{\smonth{2} \syear{2015}}% Updated by VTEXPTS2LaTeX.exe,
%27.08.2015 13:08

% ABSTRACT
%
\begin{abstract}
The paper considers functional linear regression, where scalar
responses $Y_1,\ldots,Y_n$ are modeled in dependence of i.i.d. random functions
$X_1,\ldots, X_n$.
We study
a generalization of the classical functional linear regression model.
It is assumed that there exists an unknown number of ``points of
impact,'' that is, discrete observation times where the corresponding
functional values possess significant influences on the response
variable. In addition to estimating
a functional slope parameter, the problem then is to determine the
number and locations of points of impact
as well as corresponding regression coefficients.
Identifiability of the generalized model is considered in detail. It is
shown that points of impact are
identifiable if the underlying process generating $X_1,\ldots,X_n$
possesses ``specific local variation.''
Examples are well-known processes like the Brownian motion, fractional
Brownian motion or the
Ornstein--Uhlenbeck process.
The paper then proposes an easily implementable method for estimating
the number and locations of points of impact. It is shown that this
number can be estimated consistently. Furthermore, rates of convergence for
location estimates, regression coefficients and the slope parameter are derived.
Finally, some simulation results as well as a real data application are
presented.
\end{abstract}

% KEYWORDS

%
%\begin{keyword}[class=AMS]
%\kwd[Primary ]{}
%\kwd{}
%\kwd[; secondary ]{}
%\end{keyword}
%\begin{keyword}
%\kwd{}
%\end{keyword}
\end{frontmatter}

%%%%%%%%%%%%%%%%%%%%%%%%%%%%%
%s1 #&#
\section{Introduction}\label{sec1}
%%%%%%%%%%%%%%%%%%%%%%%%%%%%%

We consider linear regression involving a scalar response variable $Y$
and a functional predictor
variable $X\in L^2([a,b])$, where $[a,b]$ is a bounded interval of
$\mathbb{R}$. It is assumed that data consist of an i.i.d. sample
$(X_i,Y_i)$, $i=1,\ldots,n$, from $(X,Y)$. The functional variable $X$
is such that $\mathbb{E}( \int_a^bX^2(t)\,dt)<+\infty$ and for
simplicity the variables are supposed to be centered in the following:
$\mathbb{E}(Y)=0$ and $\mathbb{E}(X(t))=0$ for $t\in[a,b]$ a.e.

In this paper, we study the following \textit{functional linear
regression model with points of impact}
%{\small
%e1.1 #&#
\begin{equation}
\label{impact-model} Y_i= \int_a^b
\beta(t)X_i(t)\,dt + \sum_{r=1}^S
\beta_r X_i(\tau _r)+\varepsilon_i,\qquad
i=1,\ldots,n,
\end{equation}
%
%}
where $\varepsilon_i$, $i=1,\ldots,n$ are i.i.d. centered real random
variables with
$\mathbb{E}(\varepsilon_i^2)=\sigma^2<\infty$, which are independent of
$X_i(t)$ for all $t$,
$\beta\in L^2([a,b])$ is an unknown, bounded slope function and
$ \int_a^b \beta(t) X_i(t) \,dt$ describes a common effect of the
whole trajectory $X_i(\cdot)$ on $Y_i$. In addition, the model
incorporates an
unknown number $S\in\mathbb{N}$ of ``points of impact,'' that is, \emph{specific} time points $\tau_1,\ldots,\tau_S$ with the
property that the corresponding functional values
$X_i(\tau_1),\ldots,X_i(\tau_S)$ possess some significant influence on
the response variable
$Y_i$. The function $\beta(t)$, the number $S\geq0$, as well as $\tau
_r$ and $\beta_r$,
$r=1,\ldots,S$, are unknown and have to be estimated from the data.
Throughout the paper, we will assume that all points of impact are in the
interior of the interval, $\tau_r\in(a,b)$, $r=1,\ldots,S$. Standard
functional linear regression with $S=0$ as well as the point impact model
of \citet{McKeagueSen2010}, which assumes $\beta(t)\equiv0$ and
$S=1$, are special cases
of the above model.

If $S=0$, then (\ref{impact-model}) reduces to $Y_i= \int_a^b\beta
(t)X_i(t)\,dt +\varepsilon_i$. This
model has been studied in depth in theoretical and applied statistical
literature. The most
frequently used approach for estimating $\beta(t)$ then is based on functional
principal components regression [see, e.g., \citet{FrankFriedman1993},
\citet{Bosq2000}, Cardot, Ferraty and Sarda (\citeyear{CardotFerratySarda1999}), Cardot, Mas and Sarda (\citeyear{CardotMasSarda2007}) or M\"{u}ller and Stadtm\"{u}ller
(\citeyear{MuellerStadtmueller2005}) in the context of generalized
linear models]. Rates of convergence of the estimates are derived in
\citet{HallHorowitz2007} and \citet{CaiHall2006}. Alternative approaches and
further theoretical
results can, for example, be found in Crambes, Kneip and Sarda (\citeyear{CrambesKneipSarda2009}), \citet{CardotJohannes2010}, \citet{ComteJohannes2012} or \citet{DelaigleHall2012}.

There are many successful applications of the standard linear
functional regression
model. At the same time, results are often difficult to analyze from
the points of
view of model building and substantial interpretation. The underlying
problem is that
$\int_a^b\beta(t)X_i(t)\,dt$ is a weighted average of the whole
trajectory $X_i(\cdot)$
which makes it difficult to assess specific effects of local
characteristics of the
process. This lead James, Wang and Zhu (\citeyear{JamesWangZhu2009}) to consider ``interpretable
functional regression'' by assuming that $\beta(t)=0$ for most points
$t\in[a,b]$ and
identifying subintervals of $[a,b]$ with nonzero $\beta(t)$.

A different approach based on impact points is proposed by Ferraty, Hall and Vieu (\citeyear{FerratyHallVieu2010}). For a pre-specified $q\in\mathbb{N,}$ they aim to identify
a function $g$ as well as
those design points $\tau_1,\ldots,\tau_q$ which are
``most influential'' in the sense that $g(X_i(\tau_1),\ldots,X_i(\tau
_q))$ provides a best possible prediction of $Y_i$. Nonparametric
smoothing methods are used to estimate
$g$, while $\tau_1,\ldots,\tau_q$ are selected by a cross-validation
procedure. The method is applied to data from spectroscopy, where it is
of practical interest to know
which values $X_i(t)$ have greatest influence on $Y_i$.

To our knowledge, \citet{McKeagueSen2010} are the first to explicitly
study identifiability and
estimation of a point of impact in a functional regression model. For
centered variables, their model takes the form
$Y_i=\beta X_i(\tau)+\varepsilon_i$ with a single point of impact $\tau\in
[a,b]$. The underlying process $X$ is assumed to be a fractional
Brownian motion with Hurst parameter $H$. The approach is motivated by
the analysis of
gene expression data, where a key problem is to identify individual
genes associated
with the clinical outcome. \citet{McKeagueSen2010} show that consistent
estimators are obtained by
least squares,
and that the estimator of $\tau$ has the rate of convergence $n^{-{1}/{(2H)}}$. The
coefficient $\beta$ can be estimated with a parametric rate of
convergence $n^{-{1}/{2}}$.

There also exists a link between our approach and the work of \citet{HsingRen2009} who for a given grid $t_1,\ldots,t_p$ of
observation points propose a procedure for estimating linear
combinations $m(X_i)=\sum_{j=1}^p c_j X_i(t_j)$ influencing $Y_i$.
Their approach is based on an RKHS formulation of the inverse
regression dimension-reduction problem which for any $k=1,2,3,\ldots$
allows to determine a suitable
element $(\hat c_1,\ldots,\hat c_p)^T$ of the eigenspace spanned by
the eigenvectors of the $k$ leading eigenvalues
of the empirical covariance matrix of $(X_{i}(t_1),\ldots
,X_{i}(t_p))^T$. They then show consistency of the resulting estimators
$\hat m(X_i)$ as $n,p\rightarrow\infty$ and then $k\rightarrow\infty$.
Note that (\ref{impact-model}) necessarily implies that
$Y_i=m(X_i)+\varepsilon_i$, where as $p\rightarrow\infty$ $m(X_i)$ may be
written as a linear combination as considered by \citet{HsingRen2009}.
Their method
therefore
offers a way to determine consistent estimators $\hat m(X_i)$ of
$m(X_i)$, although the structure of the estimator will not allow a
straightforward identification of model components.

Assuming a linear relationship between $Y$ and $X$, (\ref
{impact-model}) constitutes
a unified approach which incorporates the standard linear regression
model as well
as specific effects of possible point of impacts. The latter may be of
substantial
interest in many applications.

Although in this paper we concentrate on the case of unknown points of
impact, we want to emphasize that in practice also models with
pre-specified points of impact may be of potential importance. This in
particular applies to situations with a functional response variable
${\cal Y}_i(t)$, defined over the same time period $t\in[a,b]$ as
$X_i$. For a specified time
point $\tau\in[a,b]$, the standard approach [see, e.g., He, M{\"{u}}ller and Wang (\citeyear{HeMuellerWang2000})] will then assume that $Y_i:={\cal Y}_i(\tau) =\int_a^b\beta_\tau
(t)X_i(t)\,dt+\varepsilon_i$, where $\beta_\tau\in L^2([a,b])$ may vary
with $\tau$. But the value
$X_i(\tau)$ of $X_i$ at the point $\tau$ of interest may have a
specific influence, and the alternative model
$Y_i:={\cal Y}_i(\tau) =\int_a^b\beta_\tau(t)X_i(t)\,dt+\beta_1X_i(\tau
)+\varepsilon_i$ with $S=1$ and a fixed point of impact may be seen as
a promising alternative. The estimation procedure proposed in Section~\ref{sec5}
can also be applied in this situation, and theoretical results imply
that under mild conditions $\beta_1$ as well as $\beta_\tau(t)$ can be
consistently estimated with nonparametric rates of convergence.
A similar modification may be
applied in the related context of functional autoregression, where
$X_1,\ldots,X_n$
denote a stationary time series of random function, and ${\cal Y}(\tau
)\equiv X_i(\tau)$ is to be predicted from $X_{i-1}$ [see, e.g., \citet{Bosq2000}].

The focus of our work lies on developing conditions ensuring
identifiability of the
components of model (\ref{impact-model}) as well as on determining
procedures for estimating number and locations of points of impact,
regression coefficients and slope parameter.

The problem of identifiability is studied in detail in Section~\ref{sec2}. The
key assumption
is that the process possesses ``specific local variation.''
Intuitively, this means that
at least some part of the local variation of $X(t)$ in a
small neighborhood $[\tau-\epsilon, \tau+\epsilon]$ of a point $\tau\in
[a,b]$ is essentially uncorrelated with the remainder of the
trajectories outside the interval $[\tau-\epsilon, \tau+\epsilon]$.
Model (\ref{impact-model}) is uniquely identified for
all processes exhibiting specific local variation. It is also shown
that the condition
of specific local variation is surprisingly weak and only requires some suitable
approximation properties of the corresponding Karhunen--Lo\`eve basis.

Identifiability of (\ref{impact-model}) does not impose any restriction
on the degree of smoothness of
the random functions $X_i$ or of the underlying covariance function.
The same is true for
the theoretical results of Section~\ref{sec5} which yield rates of convergence
of coefficient estimates, provided that points of impact are known or
that locations can be estimated with sufficient accuracy.

But nonsmooth trajectories are advantageous when trying to identify
points of impact. In order
to define a procedure for estimating number and locations of points of impact,
we therefore restrict attention to processes whose covariance
function is nonsmooth at the diagonal. It is proved in Section~\ref{sec3} that
any such process has specific local variation. Prominent examples are
the fractional Brownian motion or the Ornstein--Uhlenbeck process. From
a practical point of view, the setting of processes with nonsmooth
trajectories covers a wide range of applications. Examples are given in
Section~\ref{sec7} and in the supplementary material [Kneip, Poss and Sarda (\citeyear{KneipPossSarda2015})],
where the methodology is applied to temperature curves and near
infrared data.

An easily implementable and computationally efficient algorithm for
estimating number and locations of points of impact is presented in
Section~\ref{SEC:4}. The basic idea is to
perform a decorrelation. Instead of regressing on $X_i(t)$, we analyze
the empirical correlation
between $Y _i$ and a process $Z_{\delta,i}(t):=X_i(t)-\frac
{1}{2}(X_i(t-\delta)+
X_i(t+\delta))$ for some $\delta>0$. For the class of processes defined
in Section~\ref{sec3},
$Z_{\delta,i}(t)$ is highly correlated with $X_i(t)$ but only possesses
extremely
weak correlations with $X_i(s)$ if $|t-s|$ is large. This implies that under
model (\ref{impact-model}) local maxima
$\widehat{\tau}_r$ of the empirical correlation between $Y_i$ and
$Z_{\delta,i}(t)$
should be found at locations close to existing points of impact. The
number $S$ is
then estimated by a cut-off criterion. It is proved that the resulting estimator
$\widehat S$ of $S$ is consistent, and we derive rates of convergence
for the estimators
$\widehat{\tau}_r$. In the special case of a fractional Brownian
motion and $S=1$,
we retrieve the basic results of \citet{McKeagueSen2010}.

In Section~\ref{sec5}, we introduce least squares estimates of $\beta(t)$ and
$\beta_r$, $r=1,\ldots, S$, based on a Karhunen--Lo\`eve decomposition.
Rates of convergence for these estimates are then derived. A simulation
study is performed in Section~\ref{sec6}, while applications to a dataset is
presented in Section~\ref{sec7}. The \hyperref[app]{Appendix} is devoted to the proofs of some of
the main results. The remaining proofs as well as the application of
our method to a second dataset are gathered in the supplementary material.

%%%%%%%%%%%%%%%%%%%%%%%%%%%%%
%s2 #&#
\section{Identifiability}\label{sec2}
%%%%%%%%%%%%%%%%%%%%%%%%%%%%%

Our setup implies that $X_1,\ldots,X_n$ are i.i.d. random functions with
the same distribution as
a generic $X\in L^2([a,b])$.
In the following, we will additionally assume that $X$ possesses a
continuous covariance function
$\sigma(t,s)$, $t,s\in[a,b]$.

In a natural way, the components of model (\ref{impact-model}) possess
different interpretations. The linear
functional $ \int_a^b \beta(t) X_i(t) \,dt$ describes a \textit{common
effect} of the
whole trajectory $X_i(\cdot)$ on $Y_i$. The additional terms $ \sum_{r=1}^S \beta_r X_i(\tau_r)$
quantify \textit{specific effects} of
the functional values $X_i(\tau_1),\ldots,X_i(\tau_S)$ at the points of
impact $\tau_1,\ldots,\tau_S$.
Identifiability of an impact point $\tau_r$ quite obviously requires
that at least some part of the
local variation of $X_i(t)$ in small neighborhoods of $\tau_r$, is
uncorrelated with the remainder of the
trajectories. This idea is formalized by introducing the concept of
``specific local variation.''

\begin{definition}\label{de1}
A process $X\in L^2([a,b])$ with continuous covariance function $\sigma
(\cdot,\cdot)$ possesses
\textit{specific local variation} if for any $t\in(a,b)$ and all
sufficiently small $\epsilon>0$ there
exists a real random variable $\zeta_{\epsilon,t}(X)$ such that
with $f_{\epsilon,t}(s):=\frac{\operatorname{cov}(X(s),\zeta_{\epsilon
,t}(X))}{\operatorname{var}(\zeta_{\epsilon,t}(X))}$ the following
conditions are satisfied:
\begin{longlist}[(iii)]
\item[(i)] $0<\operatorname{var}(\zeta_{\epsilon,t}(X))<\infty$,
\item[(ii)] $f_{\epsilon,t}(t)>0$,
\item[(iii)] $|f_{\epsilon,t}(s)|\leq(1+\epsilon)f_{\epsilon,t}(t)$ for
all $s\in[a,b]$,
\item[(iv)] $|f_{\epsilon,t}(s)|\le\epsilon\cdot f_{\epsilon,t}(t)$
for all $s\in[a,b]$ with
$s\notin[t-\epsilon,t+\epsilon]$.
\end{longlist}
\end{definition}

\begin{figure}

\includegraphics{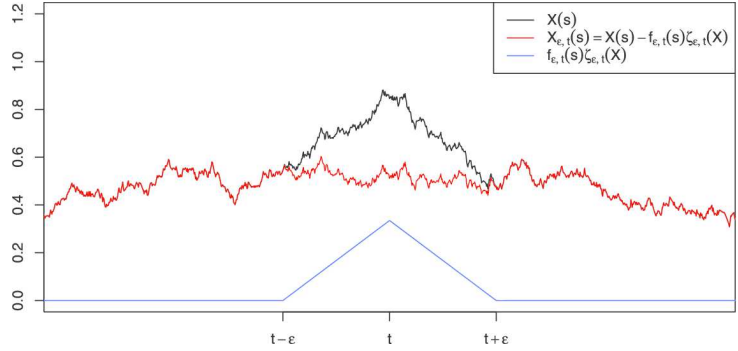}

\caption{The figure illustrates the
decomposition of a trajectory from a Brownian motion $X$ (black line)
in $X_{\epsilon,t}$ (red line) and $\zeta_{\epsilon,t}(X) f_{\epsilon
,t}$ (blue line). The component $\zeta_{\epsilon,t}(X) f_{\epsilon,t}$
can be seen to quantify the local variation of $X$ in an interval
around $t$.}
\label{fig:CovBM}
\end{figure}

 The definition of course implies that for given $t\in(a,b)$
and small $\epsilon>0$ any process $X$ with specific local variation
can be decomposed into
%
%e2.1 #&#
\begin{equation}
X(s)=X_{\epsilon,t}(s) + \zeta_{\epsilon,t}(X) f_{\epsilon,t}(s),\qquad s
\in[a,b], \label{decompeps}
\end{equation}
where $X_{\epsilon,t}(s)=X(s)-\zeta_{\epsilon,t}(X) f_{\epsilon,t}(s)$
is a process which is
uncorrelated with $\zeta_{\epsilon,t}(X)$. If $\sigma_{\epsilon,t}(\cdot
,\cdot)$ denotes the
covariance function of $X_{\epsilon,t}(s)$, then obviously
%
%e2.2 #&#
\begin{equation}
\label{decomsig} \sigma(s,u)=\sigma_{\epsilon,t}(s,u)+\operatorname{var}\bigl(
\zeta_{\epsilon,t}(X)\bigr) f_{\epsilon,t}(s) f_{\epsilon,t}(u),\qquad s,u
\in[a,b].
\end{equation}
By condition (iv), we can infer that for small $\epsilon>0$ the
component $\zeta_{\epsilon,t}(X)\times\break  f_{\epsilon,t}(s)$
essentially quantifies local variation in a small interval around the
given point $t$, since
$\frac{f_{\epsilon,t}(s)^2}{f_{\epsilon,t}(t)^2}\leq\epsilon^2$ for
all $s\notin[t-\epsilon,t+\epsilon]$.
When $X$ is a standard Brownian motion it is easily verified that
conditions (i)--(iv) are satisfied for $\zeta_{\epsilon,t}(X)=
X(t)-\frac{1}{2}(X(t-\epsilon)+X(t+\epsilon))$. Then $f_{\epsilon
,t}(s):=\frac{\operatorname{cov}(X(s),\zeta_{\epsilon,t}(X))}{\operatorname{var}(\zeta_{\epsilon,t}(X))}=1$
for $t=s$, while $f_{\epsilon,t}(s)=0$ for all $s\in[a,b]$ with
$|t-s|\geq\epsilon$. Figure~\ref{fig:CovBM} illustrates the
decomposition of $X(s)$ in $X_{\epsilon,t}(s)$ and $\zeta_{\epsilon
,t}(X) f_{\epsilon,t}(s)$ for a trajectory of a Brownian motion.
%
%f1 #&#

The following theorem shows that under our setup all impact points in
model~(\ref{impact-model}) are uniquely identified for any process
possessing specific local variation. Recall that (\ref{impact-model})
implies that
\[
m(X):=\mathbb{E}(Y|X) = \int_a^b
\beta(t)X(t)\,dt + \sum_{r=1}^S\beta
_rX(\tau_r).
\]

%th1 #&#
\begin{thmm} \label{thmident}
Under our setup, assume that $X$ possesses specific local variation.
Then, for any bounded function $\beta^*\in L^2([a,b])$, all $S^*\geq
S$, all $\beta_1^*,\ldots,\beta_{S^*}^*\in\mathbb{R}$, and
all $\tau_1,\ldots,\tau_{S^*}\in(a,b)$ with $\tau_{k}\notin\{\tau_1,\ldots
,\tau_S\}$, $k=S+1,\ldots,S^*$, we obtain
%
%e2.3 #&#
\begin{eqnarray}
\mathbb{E} \Biggl( \Biggl(m(X)-\int_a^b
\beta^*(t)X(t)\,dt - \sum_{r=1}^{S^*}
\beta_r^*X_i(\tau_r) \Biggr)^2
\Biggr)>0, \label{eqident}
\end{eqnarray}
whenever $\mathbb{E}((\int_a^b(\beta(t)-\beta^*(t))X(t)\,dt)^2)>0$,
or $\sup_{r=1,\ldots,S}|\beta_r-\beta^*_r|>0$,
or $\sup_{r=S+1,\ldots,S^*}|\beta^*_r|>0$.
\end{thmm}

The question arises whether it is possible to find general conditions
which ensure that a process
possesses specific variation. From a theoretical point of view, the
Karhunen--Lo\`eve decomposition
provides a tool for analyzing this problem.

For $f,g\in L^2([a,b])$ let $\langle f,g \rangle=
\int_a^b f(t)g(t)\,dt$ and $\|f\|$ the associated norm. We will use
$\lambda_1\geq\lambda_2\geq\cdots$ to
denote the nonzero eigenvalues of the covariance operator $\Gamma$ of
$X$, while $\psi_1,\psi_2,\ldots$
denote a corresponding system of orthonormal eigenfunctions.
It is then well known that $X$ can be decomposed
in the form
%
%e2.4 #&#
\begin{equation}
X(t)=\sum_{r=1}^\infty\langle X,
\psi_r \rangle\psi_r(t), \label{karlov}
\end{equation}
where $\mathbb{E}(\langle X,\psi_r\rangle^2)=\lambda_r$, and $\langle
X,\psi_r\rangle$ is uncorrelated with
$\langle X,\psi_l\rangle$ for $l\neq r$.

The existence of specific local variation requires that the structure
of the process is not too simple in the sense that the realizations
$X_i$ a.s. lie in a finite dimensional subspace of $L^2([a,b])$. Indeed,
if $\Gamma$ only possesses a finite number $K<\infty$ of nonzero
eigenvalues, then model (\ref{impact-model}) is not identifiable. This
is easily verified:
$X(t)=\sum_{r=1}^K \langle X,\psi_r \rangle\psi_r(t)$ implies that\vadjust{\goodbreak}
$\int_a^b\beta(t)X(t)\,dt=\sum_{r=1}^K \alpha_r \langle X,\psi_r \rangle$
with $\alpha_r=\langle\psi_r,\beta\rangle$.
Hence, there are infinitely many different
collections of $K$ points $\tau_1,\ldots,\tau_K$ and corresponding
coefficients $\beta_1,\ldots,\beta_K$ such that
\[
\int_a^b\beta(t)X(t)\,dt=\sum
_{s=1}^K \alpha_s \langle X,
\psi_s\rangle= \sum_{s=1}^K
\langle X,\psi_s \rangle \sum_{r=1}^K
\beta_r \psi_s(\tau_r)= \sum
_{r=1}^K \beta_r X(
\tau_r).
\]

Most work in functional data analysis, however, relies on the
assumption that $\Gamma$ possesses infinitely many
nonzero eigenvalues. In theoretically oriented papers, it is often
assumed that $\psi_1,\psi_2,\ldots$
form a complete orthonormal system of $L^2([a,b])$ such that $\Vert\sum_{r=1}^\infty\langle f,\psi_r \rangle
\psi_r-f\Vert=0$ for any function $f\in L^2([a,b])$.

 The following theorem shows that
$X$ possesses specific local variation if for a suitable class of functions
$L^2$-convergence generalizes to $L^\infty$-convergence.

For $t\in(a,b)$ and $\epsilon>0$, let ${\cal C}(t,\epsilon,[a,b])$
denote the space of all
continuous functions $f\in L^2([a,b])$ with the properties that
$f(t)=\sup_{s\in[a,b]}f(s)=1$ and $f(s)=0$ for $s\notin[t-\epsilon
,t+\epsilon]$.
%
%th2 #&#
\begin{thmm} \label{thmkarlov}
Let $\psi_1,\psi_2,\ldots$
be a system of orthonormal eigenfunctions corresponding to the nonzero
eigenvalues of
the covariance operator $\Gamma$ of $X$. If for all $t\in(a,b)$ there
exists an $\epsilon_t>0$ such
that
%
%e2.5 #&#
\begin{eqnarray}
\label{karlovp0} \lim_{k\rightarrow\infty}\inf_{f\in{\cal C}(t,\epsilon,[a,b])}
\sup
_{s\in[a,b]} \Biggl|f(s)-\sum_{r=1}^k
\langle f,\psi_r\rangle \psi_r(s)\Biggr|=0
\nonumber
\\[-8pt]
\\[-8pt]
\eqntext{\mbox{for every }
0<\epsilon<\epsilon_t,}
\end{eqnarray}
then the process
$X$ possesses specific local variation.
\end{thmm}

The message of the theorem is that existence of specific local
variation only requires that
the underlying basis $\psi_1,\psi_2,\ldots$ possesses suitable
approximation properties.
Somewhat surprisingly, the degree of smoothness of the realized
trajectories does not play
any role.

As an example consider a standard Brownian motion defined on
$[a,b]=[0,1]$. The corresponding Karhunen--Lo\`eve decomposition
possesses eigenvalues $\lambda_r=\frac{1}{(r-0.5)^2\pi^2}$ and
eigenfunctions\vspace*{1pt} $\psi_r(t)=\sqrt{2}\sin((r-1/2)\pi t)$, $r=1,2,\ldots.$
In the Supplementary Appendix B [Kneip, Poss and Sarda (\citeyear{KneipPossSarda2015})], it is verified
that this system of orthonormal eigenfunctions satisfies (\ref
{karlovp0}). Although all eigenfunctions are smooth,
it is well known that realized trajectories of a Brownian motion are
a.s. not
differentiable. This can be seen as a consequence of the fact that the
eigenvalues
$\lambda_r\sim\frac{1}{r^2}$ decrease fairly slowly, and, therefore,
the sequence
$\mathbb{E}( (\sum_{r=1}^k \langle X,\psi_r \rangle\psi_r'(t))^2)= \sum_{r=1}^k \lambda_r (\psi_r'(t))^2$ diverges as $k\rightarrow\infty$. At
the same time, another process with the same system of
eigenfunctions but exponentially decreasing eigenvalues $\lambda
_r^*\sim\exp(-r)$ will a.s.\vadjust{\goodbreak} show sample paths possessing an infinite
number of derivatives. Theorem~\ref{thmkarlov} states
that any process of this type still has specific local variation.

%%%%%%%%%%%%%%%%%%%%%%%%%%%%%
%s3 #&#
\section{Covariance functions which are nonsmooth at the diagonal}\label{sec3}
%%%%%%%%%%%%%%%%%%%%%%%%%%%%%

In the following, we will concentrate on developing a theoretical
framework which allows to define
an efficient procedure for estimating number and locations of points of impact.

Although specific local variation may well be present for processes
possessing very smooth sample paths, it is clear that detection of
points of impact will profit from a high local
variability which goes along with nonsmoothness. As pointed out in the
\hyperref[sec1]{Introduction}, we also believe that assuming nonsmooth trajectories
reflect the situation encountered in a number of important
applications. \citet{McKeagueSen2010} convincingly demonstrate that
genomics data lead to sample paths with fractal behavior. All important
processes analyzed in economics exhibit strong random fluctuations.
Observed temperatures
or precipitation rates show wiggly trajectories over time, as can be
seen in our application in Section~\ref{sec7}. Furthermore, any growth process
will to some extent be influenced by random changes in environmental
conditions. In functional data analysis, it is common practice to
smooth observed (discrete) sample paths and to interpret nonsmooth
components as
``errors.'' We want to emphasize that, unless observations are
inaccurate and there exists some important measurement error, such
components are an intrinsic part of the process. For many purposes as,
for example, functional principal component analysis, smoothing makes a
lot of sense since local variation
has to be seen as nuisance. But in the present context local variation
actually is a key property for identifying impact points.

Therefore, further development will focus on processes with nonsmooth
sample paths which will be
expressed in terms of a nonsmooth diagonal of the corresponding
covariance function $\sigma(t,s)$.
It will be assumed that $\sigma(t,s)$ possesses nonsmooth trajectories
when passing from $\sigma(t,t-\Delta)$ to $\sigma(t,t+\Delta)$, but is twice
continuously differentiable for all $(t,s)$, $t\neq s$. An example is
the standard Brownian motion
whose covariance function $\sigma(t,s)=\min(t,s)$ has a kink at the
diagonal. Indeed, in view
of decomposition (\ref{decomsig}) a nonsmooth transition at diagonal
may be seen as a natural consequence of
pronounced specific local variation.

For a precise analysis, it will be useful to reparametrize the
covariance function. Obviously, the
symmetry of $\sigma(t,s)$ implies that
\begin{eqnarray*}
\sigma(t,s)&=&\sigma\bigl(\tfrac{1}{2}\bigl(t+s+|t-s|\bigr),\tfrac
{1}{2}\bigl(t+s-|t-s|\bigr)
\bigr)\\
&=:&\omega^*\bigl(t+s,|t-s|\bigr) \qquad\mbox{for all } t,s\in[a,b].
\end{eqnarray*}
Instead of $\sigma(t,s)$, we may thus equivalently consider the
function $\omega^*(x,y)$ with $x=t+s$ and
$y=|t-s|$. When passing from $s=t-\Delta$ to $s=t+\Delta$,
the degree of smoothness of $\sigma(t,s)$ at $s=t$ is reflected by the
behavior of
$\omega^*(2t,y)$ as $y\rightarrow0$.

First, consider the case that $\sigma$ is twice continuously
differentiable and for fixed $x$ and $y>0$ let
$\frac{\partial}{\partial y_+}\omega^*(x,y)|_{y=0}$ denote the right
(partial) derivative of $\omega^*(x,y)$ as $y\rightarrow0$.
It is easy to check that in this
case for all $t\in(a,b)$ we obtain
%
%e3.1 #&#
\begin{eqnarray}\label{der1}
\frac{\partial}{\partial y_+}\omega^*(2t,y)\bigg|_{y=0} &=&\frac{\partial
}{\partial y} \sigma
\biggl(t+\frac{y}{2},t-\frac{y}{2}\biggr)\bigg|_{y=0}
\nonumber
\\[-8pt]
\\[-8pt]
\nonumber
&=&
\frac{1}{2}\biggl(\frac{\partial}{\partial s} \sigma(s,t)\bigg|_{s=t}-
\frac{\partial}{\partial s} \sigma(t,s)\bigg|_{s=t}\biggr)=0.
\end{eqnarray}

In contrast, any process with $\frac{\partial}{\partial y_+}\omega
^*(x,y)|_{y=0}\neq0$ is nonsmooth at the diagonal. If this function
is smooth for all other points $(x,y)$, $y>0$, then the process,
similar to the Brownian motion, possesses a kink at the diagonal.
Now note that, for any process with $\sigma(t,s)=\omega^*(t+s,|t-s|)$
continuously differentiable for $t\neq s$ but $\frac{\partial}{\partial
y_+}\omega^*(x,y)|_{y=0}<0$, it is possible to find a twice
continuously differentiable function $\omega(x,y,z)$ with $\sigma
(t,s)=\omega(t,s,|t-s|)$ such that
$\frac{\partial}{\partial y_+}\omega^*(t+t,y)|_{y=0}=\frac{\partial
}{\partial y}\omega(t,t,y)|_{y=0}$.

In a still more general setup, the above ideas are formalized by
Assumption~\ref{assum1}
below which, as will be shown in Theorem~\ref{thmcharX}, provides
sufficient conditions in order to guarantee that the underlying process
$X$ possesses
specific variation. We will also allow for unbounded derivatives as
$|t-s|\rightarrow0$.

%as1 #&#
\begin{assumption} \label{assum1}
For some open subset $\Omega\subset\mathbb{R}^3$ with $[a,b]^2\times
[0,b-a]\subset\Omega$,
there exists a twice continuously differentiable function $\omega
:\Omega\rightarrow
\mathbb{R}$ as well as some $0<\kappa<2$ such that for all $t,s\in[a,b]$
%
%e3.2 #&#
\begin{equation}
\sigma(t,s)=\omega\bigl(t,s,|t-s|^\kappa\bigr). \label{deromega}
\end{equation}
Moreover,
%
%e3.3 #&#
\begin{equation}
0<\inf_{t\in[a,b]}c(t)\qquad \mbox{where } c(t):=-\frac{\partial}{\partial z}
\omega(t,t,z)\bigg|_{z=0}. \label{deromega2}
\end{equation}
\end{assumption}

One can infer from (\ref{der1}) that for every twice continuously
differentiable covariance function $\sigma$
there exists some function $\omega$ such that (\ref{deromega}) holds
with $\kappa=2$.
But note that formally introducing $|t-s|^\kappa$ as an extra argument
establishes an easy way of capturing nonsmooth behavior as
$|t-s|\rightarrow0$, since $\sigma$ is not twice differentiable at the
diagonal if $\kappa<2$. In
Assumption~\ref{assum1}, the value of $\kappa<2$ thus quantifies the
degree of smoothness of
$\sigma$ at the diagonal. A very small $\kappa$ will reflect pronounced
local variability and
extremely nonsmooth sample paths.
There are many well-known processes satisfying this assumption.

%\begin{itemize}
%\item
\textit{Fractional Brownian motion} with Hurst coefficient $0<H<1$ on an
interval $[a,b]$, $a> 0$: The covariance function is then
given by
\[
\sigma(t,s)=\tfrac{1}{2}\bigl(t^{2H}+s^{2H}-|t-s|^{2H}
\bigr).
\]
In this case, Assumption~\ref{assum1} is satisfied with $\kappa=2H$,
$\omega(t,s,z)=
\frac{1}{2}(t^{2H}+s^{2H}-z)$
and $c(t)=1/2$.

%\item
\textit{Ornstein--Uhlenbeck process} with parameters $\sigma_{u}^2,\theta
>0$: The covariance function is then defined by
\[
\sigma(t,s)=\frac{\sigma_{u}^2}{2\theta}\bigl(\exp\bigl(-\theta|t-s|\bigr)-\exp\bigl(-\theta(t+s)\bigr)\bigr).
\]
Then Assumption~\ref{assum1} is satisfied with $\kappa=1$,
$\omega(t,s,z)=
\frac{\sigma_{u}^2}{2\theta}(\exp(-\theta z)-\break  \exp(-\theta(t+s)))$
and $c(t)=\sigma_{u}^2/2$.
%\end{itemize}

Theorem~\ref{thmcharX} below now states that any process respecting
Assumption~\ref{assum1} possesses specific local variation. In Section~\ref{sec2}, we already discussed the structure of
an appropriate r.v. $\zeta_{\epsilon,t}(X)$ for the special case of a
standard Brownian
motion. The same type of functional may now be used in a more general setting.

For $\delta>0$ and $[t-\delta,t+\delta]\subset[a,b]$, define
%
%e3.4 #&#
\begin{equation}
\label{Zdef} Z_{\delta}(X,t) = X(t) - \tfrac{1}{2} \bigl(X(t-
\delta)+X(t+\delta ) \bigr).
\end{equation}

%th3 #&#
\begin{thmm} \label{thmcharX}
Under our setup, assume that the covariance function $\sigma$ of $X$
satisfies Assumption~\ref{assum1}. Then
$X$ possesses specific local variation, and for any $\epsilon>0$ there
exists a $\delta>0$
such that conditions \textup{(i)--(iv)} of Definition~\ref{de1} are satisfied for
$\zeta_{\epsilon,t}(X)=Z_{\delta}(X,t)$, where $Z_{\delta}(X,t)$ is
defined by
(\ref{Zdef}).
\end{thmm}

%%%%%%%%%%%%%%%%%%%%%%%%%%%%%
%s4 #&#
\section{Estimating points of impact}\label{SEC:4}
%%%%%%%%%%%%%%%%%%%%%%%%%%%%%

When analyzing model (\ref{impact-model}), a central problem is to
estimate number and locations of
points of impact.
Recall that we assume an i.i.d. sample $(X_i,Y_i)$, $i=1,\ldots,n$,
where $X_i$ possesses the
same distribution as a generic $X$. Furthermore, we consider the case
that each $X_i$ is
evaluated at $p$ equidistant points $t_j=a+\frac{j-1}{p-1}(b-a)$, $j=1,
\ldots, p$.

\begin{remark*}
Note that all variables have been assumed to have means
equal to zero. Any practical application of the methodology introduced below,
however, should rely on centered data to be obtained from the original
data by subtracting sample means. Obviously, the theoretical results
developed in this section remain unchanged for this situation with
however substantially longer proofs.
\end{remark*}

Determining $\tau_1,\ldots,\tau_S$ of course constitutes a model
selection problem. Since in practice
the random functions $X_i$ are observed on a discretized grid of $p$
points, one may tend to use
multivariate model selection procedures like Lasso or related methods.
But these procedures
are multivariate in nature and are not well adapted to a functional
context. An obvious difficulty is the linear functional
$ \int_a^b \beta(t) X_i(t) \,dt\approx\frac{1}{p} \sum_{j=1}^p
\beta(t_j) X_i(t_j)$ which contradicts
the usual sparseness assumption by introducing some common effects of
all variables.
But even if $ \int_a^b \beta(t) X_i(t) \,dt\equiv0$,
results may heavily depend on the number $p$ of observations per function.
Note that in our functional setup for any fixed $m\in\mathbb{N}$ we
necessarily have $\operatorname{Var}(X_i(t_j)-X_i(t_{j-m}))\rightarrow
0$ as $p\rightarrow\infty$.
Lasso theory, however, is based on the assumption that variables are
not too heavily correlated. For example, the results of Bickel, Ritov and Tsybakov
(\citeyear{Bickel2009}) indicate
that convergence of parameter estimates
\textit{at least} requires that $\sqrt{n/\log
p}(\operatorname{Var}(X_i(t_j)-X_i(t_{j-1})))\rightarrow\infty$ as $n\rightarrow
\infty$. This follows from the distribution version of the restricted
eigenvalue assumption and Theorem~5.2 of Bickel, Ritov and Tsybakov (\citeyear{Bickel2009}) [see also
Zhou, van~de Geer and B\"uhlmann (\citeyear{ZhouvandeGeerBuehlmann2009}) for a discussion on correlation assumptions for
selection models]. As a consequence, standard multivariate model
selection procedures cannot work unless the number
$p$ of grid points is sufficiently small compared to $n$.

In this paper, we propose a very simple approach which is based on the
concepts developed in the preceding sections. The idea is to identify
points of impact by determining the grid points
$t_j$, where $Z_{\delta,i}(t_j):=Z_\delta(X_i,t_j)$ possesses a
particularly high correlation with
$Y_i$.

The motivation of this approach is easily seen when considering our
regression model (\ref{impact-model}) more closely. Note that $Z_{\delta
,i}(t)$ is strongly
correlated with $X_i(t)$, but it is ``almost'' uncorrelated with
$X_i(s)$ for
$|t-s|\gg\delta$. This in turn implies that the correlation between
$Y_i$ and
$Z_{\delta,i}(t)$ will be comparably high if and only if a particular
point $t$ is
close to a point of impact. More precisely, Lemmas 3 and 4 in the
Supplementary Appendix C [Kneip, Poss and Sarda (\citeyear{KneipPossSarda2015})] show
that
as $\delta\rightarrow0$ and $\min_{r\neq s} |\tau_s-\tau_r|\gg\delta$
\begin{eqnarray*}
\mathbb{E} \bigl(Z_{\delta,i}(t_j)Y_i \bigr)&=&
\beta_r c(\tau_r)\delta ^\kappa+O\bigl(\max\bigl
\{\delta^{\kappa+1},\delta^2\bigr\}\bigr)\qquad \mbox{if }
|t_j-\tau_r|\approx0,
\\
\mathbb{E} \bigl(Z_{\delta,i}(t_j)Y_i \bigr)&=&O
\bigl(\max\bigl\{\delta^{\kappa
+1},\delta^2\bigr\}\bigr)\qquad \mbox{if } \min_{r=1,\ldots,S} |t_j-\tau _r| \gg\delta.
\end{eqnarray*}
Moreover, assuming that the process $X$ possesses a Gaussian
distribution, then since $\operatorname{Var}(Z_{\delta,i}(t_j))=O(\delta^\kappa)$ [see
(\ref{Zdelta2}) in the proof of Theorem~\ref{thmcharX}], the
Cauchy--Schwarz inequality lead to
$\operatorname{Var}(Z_{\delta,i}(t_j)Y_i)=O(\delta^\kappa)$, and hence
\[
\Biggl|\frac{1}{n}\sum_{i=1}^n
Z_{\delta,i}(t_j)Y_i-\mathbb{E}
\bigl(Z_{\delta
,i}(t_j)Y_i\bigr)\Biggr|=O_P
\biggl(\sqrt{\frac{\delta^\kappa}{n}}\biggr).
\]
These arguments indicate that points of impact may be estimated by
using the locations of sufficiently large local maxima of $|\frac
{1}{n}\sum_{i=1}^n Z_{\delta,i}(t_j)Y_i|$.
A sensible identification will require a suitable choice of $\delta>0$
in dependence
of the sample size~$n$. If $\delta$ is too large, it will not be
possible to distinguish
between the influence of points of impact which are close to each
other. On the other
hand, if $\delta$ is too small compared to $n$ (as, e.g., $\delta^k\sim
n^{-1}$), then ``true'' maxima may perish in a flood of random peaks.

\begin{figure}

\includegraphics{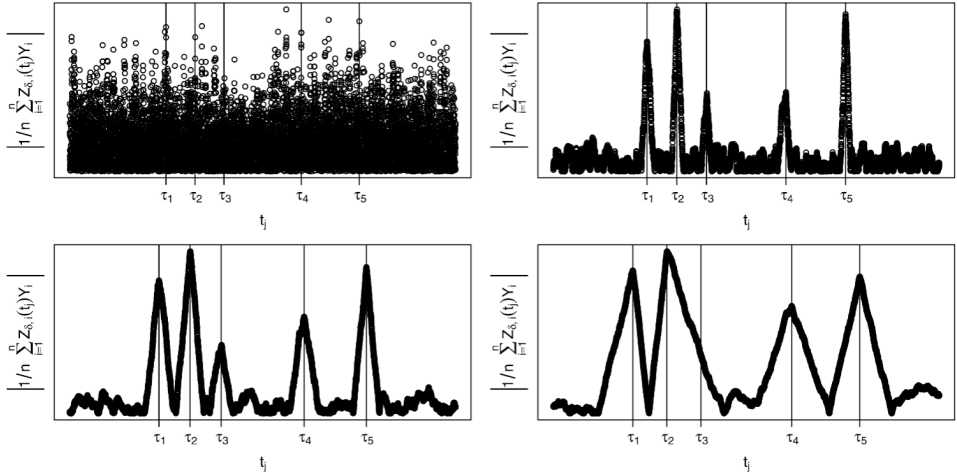}

\caption{The figure shows $|\frac
{1}{n}\sum_{i=1}^n Z_{\delta,i}(t_j)Y_i|$ for different
choices of $\delta$ in a point of impact model with $5$ points of
impact whose locations are indicated by vertical lines. The upper left
panel corresponds to a very small $\delta$, where the noise level
overlays the signal. By increasing $\delta$ the location of the points
of impact becomes more and more visible. % as the absolute correlation
%around the point of impacts increases.
By choosing $\delta$ too large, as in the lower right panel, we are not
able to distinguish between the influence of points of impact in close
vicinity anymore.}
\label{fig:1}
\end{figure}

The situation is illustrated in Figure~\ref{fig:1}. It shows a
simulated example
of the regression model (\ref{impact-model}) with $n=5000$, $\beta
(t)\equiv0$, and $S=5$ points
of impact. The error term is standard normal, while $X_i$ are
independent realizations
of an Ornstein--Uhlenbeck process with $\theta= 5$ and $\sigma_u =
3.5$, evaluated over $p=10\mbox{,}001$ equidistant grid points in the interval
$[0,1]$. The figure shows the behavior of $|\frac{1}{n}\sum_{i=1}^n
Z_{\delta,i}(t_j)Y_i|$ for different choices
$\delta=10/10\mbox{,}001\approx5/n$,
$\delta=142/10\mbox{,}001\approx1/\sqrt{n}$,
$\delta=350/10\mbox{,}001\approx2.47/\sqrt{n}$, and
$\delta=750/10\mbox{,}001\approx5.3/\sqrt{n}$.

In order to consistently estimate $S$,
our estimation procedure requires to exclude all points $t$ in an
interval of size $\sqrt{\delta}$ around the local maxima of $|\frac
{1}{n}\sum_{i=1}^n Z_{\delta,i}(t_j)Y_i|$ from further considerations.
The vertical lines in Figure~\ref{fig:1} indicate the true location of
the points of impact, whereas the tick marks on the horizontal axis
represent our possible candidates for $\tau$ when applying the
following estimation procedure.
%
%f2 #&#

\textit{Estimation procedure}:
Choose some $\delta>0$ such that there exists some $k_\delta\in\mathbb
{N}$ with $1\leq k_\delta<\frac{p-1}{2}$
and $\delta=k_\delta(b-a)/(p-1)$.
In a first step, determine for all $j\in{\cal J}_{0,\delta}:=\{
k_\delta+1,\ldots,p-k_\delta\}$
\[
Z_{\delta,i}(t_j):=X_i(t_j)-
\tfrac{1}{2}\bigl(X_i(t_j-\delta)+X_i(t_j+
\delta)\bigr).
\]
Iterate for $l=1,2,3,\ldots:$
\begin{itemize}
\item Determine
\[
j_l=\arg \max_{j\in{\cal J}_{l-1,\delta}} \Biggl|\frac{1}{n}\sum
_{i=1}^n Z_{\delta,i}(t_j)Y_i\Biggr|
\]
and set $\widehat{\tau}_l:=t_{j_l}$.
\item Set ${\cal J}_{l,\delta}:=\{j\in{\cal J}_{l-1,\delta}|
|t_j-\widehat{\tau}_l|\ge\sqrt{\delta}/2\}$, that is, eliminate all
points in an interval
of size $\sqrt{\delta}$ around $\widehat{\tau}_l$. Stop iteration if
${\cal J}_{l,\delta}=\varnothing$.
\end{itemize}
Choose a suitable cut-off parameter $\lambda>0$.
\begin{itemize}
\item Estimate $S$ by
\[
\widehat S =\arg \min_{l=0,1,2,\ldots} \biggl|\frac{({1}/{n})\sum_{i=1}^n
Z_{\delta,i}(\widehat{\tau}_{l+1})Y_i}{
(({1}/{n})\sum_{i=1}^n Z_{\delta,i}(\widehat{\tau
}_{l+1})^2)^{1/2}}\biggr|< \lambda.
\]
\item $\widehat{\tau}_1,\ldots,\widehat{\tau}_{\widehat S}$ then are
the final estimates of the points of impact.
\end{itemize}

A theoretical justification for this estimation procedure is given by
Theorem~\ref{poicon}.
Its proof along with the proofs of Propositions~\ref{poiconkappa} and \ref{lem4} below can be
found in the Supplementary Appendix C.
Theory relies on an asymptotics $n\rightarrow\infty$ with $p\equiv
p_n \geq
L n^{1/\kappa}$ for some constant $0<L<\infty$.
It is based on the following additional assumption
on the structure of $X$ and $Y$.

%as2 #&#
\begin{assumption} \label{assum2}
(a) $X_1,\ldots,X_n$ are i.i.d. random functions distributed
according to $X$. The process $X$ is
\textit{Gaussian} with covariance function $\sigma(t,s)$.

(b) The error terms $\varepsilon_1,\ldots,\varepsilon_n$ are i.i.d.
$N(0,\sigma^2)$ r.v. which are
independent of~$X_i$.
\end{assumption}

%th4 #&#
\begin{thmm} \label{poicon}
Under our setup and Assumptions \ref{assum1} as well as \ref{assum2}
let $\delta\equiv\delta_n\rightarrow0$ as $n\rightarrow\infty$ such that
$\frac{n\delta^\kappa}{|\log\delta|}\rightarrow\infty$ as well as
$\frac{\delta^\kappa}{n^{-\kappa+1}}\rightarrow0$. As $n\rightarrow
\infty$ we then obtain
%
%e4.1 #&#
\begin{equation}
\label{thmm3eq1} \max_{r=1,\ldots,\widehat{S}}\min_{s=1,\ldots,S}|
\widehat{\tau}_r-\tau _s| = O_P\bigl(
n^{-{1}/{k}}\bigr).
\end{equation}
Additionally, assume that
$\delta^2=O(n^{-1})$ and that the algorithm is applied with cut-off parameter
\[
\lambda\equiv\lambda_n=A\sqrt{ \frac{\operatorname{Var}(Y_i)}{n} \log\biggl(
\frac
{b-a}{\delta}\biggr)} \qquad\mbox{where } A>\sqrt{2}.
\]
Then
%
%e4.2 #&#
\begin{equation}
\label{thmm3eq4} P(\widehat{S}=S) \rightarrow 1\qquad \mbox{as } n\rightarrow\infty.
\end{equation}
\end{thmm}

The theorem of course implies that the rates of convergence of\break the
estimated points of impact depend on $\kappa$. If $\kappa=1$ as, for\break 
example, for the Brownian motion or the Ornstein--Uhlenbeck process, then\break 
$\max_{r=1,\ldots,\widehat S} \min_{s=1,\ldots,S} |\widehat{\tau}_r-\tau
_s|=O_P( n^{-1})$. Arbitrarily fast rates
of convergence can be achieved for very nonsmooth processes with
$\kappa\ll1$.

A suitable choice of $\delta$ satisfying the requirements of the theorem
for all possible $\kappa<2$ is $\delta=C n^{-1/2}$ for some constant $C$.

Recall that for $l>1$, our algorithm requires that $\widehat{\tau}_l$
is determined
only from those points $t_j$ which are not in $\sqrt{\delta
}/2$-neighborhoods of
any previously selected $\widehat{\tau}_1,\ldots,\widehat{\tau}_{l-1}$.
This implies that
for any $\delta$ the number $M_\delta$ of iteration steps is finite, and
$M_\delta=O(\frac{b-a}{\sqrt{\delta}/2})$ is the maximal possible
number of ``candidate'' impact points which can be detected for a fixed
$n$ and $\delta\equiv\delta_n$.
The size of these intervals
is due to the use of the cut-off criterion for estimating $S$.
It can easily be seen
from the proof of the theorem that in order to establish (\ref
{thmm3eq1}) it suffices
to eliminate all points in $\delta|\log\delta|$ neighborhoods of
$\widehat{\tau}_1,\ldots,\widehat{\tau}_{l-1}$ which is a much weaker
restriction.

We also want to emphasize that the cut-off value provided by the
theorem heavily relies
on the Gaussian assumption. A different approach that may work under
more general
conditions is to consider all selected local maxima $\widehat{\tau
}_1,\ldots,\widehat{\tau}_{M_\delta}$ and to estimate $S$ by usual model
selection
criteria like BIC.

This is quite easily done if it can additionally be assumed that, in
model (\ref{impact-model}),
$\beta(t)=0$ for all $t\in[a,b]$. One may then
apply a best subset selection by regressing $Y_i$ on all possible
subsets of $X_i(\widehat{\tau}_1),\ldots,X_i(\widehat{\tau}_{M_\delta
})$, and by calculating
the residual sum of squares $\mathit{RSS}_s$
for each subset of size $s$. An
estimate $\widehat S$ is obtained by minimizing
%
%e4.3 #&#
\begin{equation}
\label{biccriterion} \mathit{BIC}_s = n\log(\mathit{RSS}_s/n) + s \log(n)
\end{equation}
over all possible values of $s$.

If $\int_a^b\beta(t)X_i(t)\,dt\neq0$, this approach will of course lead
to biased
results, since part of the influence of this component on the response
variable $Y_i$
may be approximated by adding additional artificial ``points of impact.''
But
an obvious idea is then to incorporate estimates of the linear functional
by relying on functional principal components.
Recall the Karhunen--Lo\`eve decomposition already discussed in
Section~\ref{sec2},
and note that $\int_a^b\beta(t)X_i(t)\,dt =\sum_{r=1}^\infty\alpha_r
\langle X,\psi_r \rangle$ with $\alpha_r= \langle\psi_r,\beta\rangle
$. For $k,S\in\mathbb{N}$, estimates $\widehat{\psi}_r$ of
$\psi_r$ and a subset
$\tilde\tau_1,\ldots,\tilde\tau_S\in\{\widehat{\tau}_1,\ldots,\widehat
{\tau}_{M_\delta}\}$
one may consider an approximate relationship which resembles an
``augmented model'' as proposed
by \citet{KneipSarda2011} in a different context:
%
%e4.4 #&#
\begin{equation}
\label{augmentedmodel} Y_i\approx\sum_{r=1}^k
\alpha_r \langle X_i,\widehat{\psi}_r
\rangle +\sum_{r=1}^S\beta_r
X_i(\tilde\tau_r)+\varepsilon_i^*.
\end{equation}
Based on corresponding least-squares estimates of the coefficients
$\alpha_r$ and $\beta_r$,
the number $S$ and an optimal value of $k$ may then be estimated by
the BIC criterion.

This approach also offers a way to select a sensible value of $\delta
=C n^{-1/2}$ for
a suitable range of values $C\in[C_{\mathrm{min}},C_{\mathrm{max}}]$. For finite $n$,
different choices
of $C$ (and~$\delta$) may of course lead to different candidate values
$\hat\tau_r$, $r=1,2,\ldots.$ A~straightforward approach is then to
choose the value of $\delta$, where the respective estimates of impact
points lead to the best fitting augmented model (\ref{augmentedmodel}).
In addition to estimating $S$ and an optimal value of $k$, BIC
may thus also be used to approximate an optimal value of $C$ (and
$\delta$).

Recall that the above approach is applicable if Assumption~\ref{assum1}
holds for some $\kappa<2$. In a
practical application, one may thus want to check the applicability of
the theory by estimating the value
of $\kappa$ from the data. We have
$\mathbb{E}(Z_{\delta,i}(t_j)^2)=\delta^\kappa (2c(t_j)-\frac
{2^\kappa}{2}c(t_j) )+o(\delta^\kappa)$ [see (\ref{Zdelta2}) in
the proof of Theorem~\ref{thmcharX}]. Consequently, $\frac{\mathbb
{E}(Z_{\delta,i}(t_j)^2)}{\mathbb{E}(Z_{\delta/2,i}(t_j)^2)}=2^\kappa
+o(1)$ as
$\delta\rightarrow0$. Without restriction assume that $k_\delta$ is an
even number. The
above arguments motivate the estimator
\[
\label{estkappa} \widehat{\kappa}= \log_2 \biggl( \frac{({1}/{(p-2k_\delta)})\sum_{j\in{\cal J}_{0,\delta}}\sum_{i=1}^n
Z_{\delta,i}(t_j)^2}{
({1}/{(p-2k_\delta)})\sum_{j\in{\cal J}_{0,\delta}}\sum_{i=1}^n
Z_{\delta/2,i}(t_j)^2}
\biggr)
\]
of $\kappa$. In Proposition~\ref{poiconkappa} below, it is shown that
$\widehat{\kappa}$ is a consistent estimator of $\kappa$
as $n\rightarrow\infty$, $\delta\rightarrow0$. In practice, an
estimate $\widehat{\kappa}\ll2$ will
indicate a process whose covariance function possesses a nonsmooth diagonal.

%pr1 #&#
\begin{prop}\label{poiconkappa}
Under the conditions of Theorem~\ref{poicon}, we have
%
%e4.5 #&#
\begin{equation}
\label{thmm3eq3} \widehat{\kappa}=\kappa+O_P\bigl(n^{-1/2}+
\delta^{\min\{2,2/\kappa\}}\bigr).
\end{equation}
\end{prop}

A final theoretical result concerns the distance between $X_i(\widehat
{\tau}_r)$ and
$X_i(\tau_r)$. It will be of crucial importance in the next section on
parameter estimation. Without restriction, we will in the following
assume that points of impact are ordered
in such a way that $\tau_r =\arg \min_{s=1,\ldots,S} |\widehat{\tau
}_r-\tau_s|$, $r=1,\ldots,S$.

%pr2 #&#
\begin{prop}\label{lem4}
Under the assumptions of Theorem~\ref{poicon},
we obtain for every $r=1,\ldots,S$
%
%e4.6 #&#
%e4.7 #&#
\begin{eqnarray}
\frac{1}{n}\sum_{i=1}^n
\bigl(X_i(\tau_r)-X_i(\widehat{
\tau}_r)\bigr)^2&=&O_p\bigl(n^{-1}
\bigr), \label{lem4eq1}
\\
\frac{1}{n}\sum_{i=1}^n
\bigl(X_i(\tau_r)-X_i(\widehat{
\tau}_r)\bigr)\varepsilon _i&=&O_p
\bigl(n^{-1}\bigr). \label{lem4eq3}
\end{eqnarray}
\end{prop}

%%%%%%%%%%%%%%%%%%%%%%%%%%%%%
%s5 #&#
\section{Parameter estimates}\label{sec5}
%%%%%%%%%%%%%%%%%%%%%%%%%%%%%

Recall that Assumption~\ref{assum1} is only a sufficient, not a
necessary condition of identifiability.
Even if this assumption is violated and the covariance function $\sigma
(t,s)$ is very smooth, there may exist alternative procedures leading
to sensible
estimators $\widehat{\tau}_r$. In the following, we will thus only
assume that the points of impacts
are estimated by some procedure such that $P(\widehat{S}=S)\rightarrow
1$ as $n\rightarrow\infty$ and
such that (\ref{lem4eq1}) as well as (\ref{lem4eq3}) hold for all
$r=1,\ldots,S$. Note
that this assumption is trivially satisfied if analysis is based on
pre-specified points of impact as discussed in the \hyperref[sec1]{Introduction}.

In situations where it can be assumed that $\int_a^b \beta
(t)X_i(t)\,dt=0$ a.s.,
we have $Y_i=\sum_{r=1}^S \beta_r X_i(\tau_r)+
\varepsilon_i$, $i=1,\ldots,n$, and the regression
coefficient may be obtained by least squares when replacing the
unknown points of
impact $\tau_r$ by their estimates $\widehat{\tau}_r$.
More precisely, in this
case an estimator $\widehat{\bolds{\beta}}=(\widehat{\beta}_1,\ldots
,\widehat{\beta}_{\widehat{S}})^T$
of $\bolds{\beta}=(\beta_1,\ldots,\beta_S)^T$ is determined by minimizing
%
%e5.1 #&#
\begin{equation}
\label{lsbeta0} \frac{1}{n} \sum_{i=1}^n
\Biggl(Y_i-\sum_{r=1}^{\widehat{S}}
b_r X_i(\widehat {\tau}_r)
\Biggr)^2
\end{equation}
over all possible values $b_1,\ldots,b_{\widehat{S}}$.

Let $\mathbf{X}_i(\bolds{\tau}):=(X_i(\tau_1),\ldots,X_i(\tau
_S))^T$, and let
$\Sigma_\tau:=\mathbb{E}(\mathbf{X}_i(\bolds{\tau})\mathbf
{X}_i(\bolds{\tau})^T)$. Note
that identifiability of the regression model as stated in Theorem~\ref
{thmident} in
particular implies that $\Sigma_\tau$ is invertible.

If $\widehat{S}=S$, then by (\ref{lem4eq1}) and (\ref{lem4eq3}) the
differences between $\widehat{\tau}_r$
and $\tau_r$, $r=1,\ldots,S$ are asymptotically negligible, and the
asymptotic distribution of $\widehat{\bolds{\beta}}$
coincides with the asymptotic distribution the least squares estimator
to be obtained if points of impact were known:
%
%e5.2 #&#
\begin{equation}
\label{thmpNLFeq1} \sqrt{n}(\widehat{\bolds{\beta}}-\bolds{\beta})
\rightarrow_D N\bigl(0,\sigma^2 \Sigma_\tau^{-1}
\bigr)
\end{equation}
as $n\rightarrow\infty$. A proof is straightforward, and thus omitted.

In the general case with $\beta(t)\neq0$ for some $t$, we propose to
rely on the augmented model (\ref{augmentedmodel}). Thus, let $\hat
\lambda_1\geq\hat\lambda_2\geq\cdots$ and $\widehat{\psi}_1,\widehat{\psi
}_2,\ldots$ denote eigenvalues and eigenfunctions of the empirical
covariance operator of $X_1,\ldots,X_n$. Given estimates
$\widehat{\tau}_1,\ldots,\widehat{\tau}_{\widehat{S}}$ and a suitable
cut-off parameter $k$ estimates
$\widehat{\bolds{\beta}}=(\widehat{\beta}_1,\ldots,\widehat{\beta
}_{\widehat{S}})^T$
of $\bolds{\beta}=(\beta_1,\ldots,\beta_S)^T$ and $\widehat{\alpha
}_1,\ldots,\widehat{\alpha}_k$ of
$\alpha_1,\ldots,\alpha_k$ are determined by minimizing
%
%e5.3 #&#
\begin{equation}
\label{augmentedLS} \sum_{i=1}^n \Biggl(
Y_i-\sum_{r=1}^k
a_r \langle X_i,\widehat{\psi}_r \rangle-
\sum_{r=1}^{\widehat{S}}b_r
X_i(\widehat\tau_r) \Biggr)^2
\end{equation}
over all $a_r, b_s$, $r=1,\ldots,k$, $s=1,\ldots,\widehat{S}$. Based on
the estimated coefficients
$\widehat{\alpha}_1,\ldots,\widehat{\alpha}_k$, and estimator of the
slope function $\beta$ is then given
by $\widehat{\beta}(t):=\sum_{r=1}^k \widehat{\alpha}_k \widehat{\psi}_r(t)$.

In the following we will rely on a slight change of notation in the
sense that $Y_i$, $X_i$ (and $\epsilon_i$) are centered data obtained
for each case by subtracting sample means. As pointed out in the
remark, we argue that theoretical results stated in Section~\ref{SEC:4} remain
unchanged for this situation.
In the context of (\ref{augmentedLS}) centering ensures that $X_i$,
$i=1,\ldots,n$, can be \textit{exactly}
represented by $X_i=\sum_{j=1}^n\langle X_i,\widehat{\psi}_r \rangle
\widehat{\psi}_r$ (necessarily $\hat\lambda_j=0$ for $j>n$).

Our theoretical analysis of the estimators defined by (\ref
{augmentedLS}) relies on the work of \citet{HallHorowitz2007} who
derive rates of convergence of the estimator $\widehat{\beta}(t)$ in a
standard functional regression model with $S=0$. Under our Assumption~\ref{assum2} their results are additionally based on the following
assumption on the eigendecompositions of $X$ and $\beta$.
%
%as3 #&#
\begin{assumption} \label{assum3}
(a) There exist some $\mu>1$ and some $\sigma^2<C_0<\infty$ such
that $\lambda_j-\lambda_{j+1}\geq C_0^{-1} j^{-\mu-1}$ for all
$j\geq1$.

(b) $\beta(t)=\sum_{j=1}^\infty\alpha_j \psi(t)$ for all $t$,
and $|\alpha_j|\geq C_0 j^{-\nu}$ for some $\nu>1+\frac{1}{2}\mu$.
\end{assumption}
\citet{HallHorowitz2007} show that if $S=0$ and $k=O(n^{1/(\mu+2\nu
)})$, then $\int_a^b (\widehat{\beta}(t)-\beta(t))^2\,dt =O_p(n^{-(2\nu
-1)/(\mu+2\nu)})$. This is known to be an optimal rate
of convergence under the standard model.

When dealing with points of impact, some additional conditions are
required. Note that
$\sigma(t,s)=\sum_{j=1}^\infty\lambda_j \psi_j(t)\psi_j(s)$.
Let $\sigma^{[k]}(t,s):=\sum_{j=k+1}^\infty\lambda_j\times\break  \psi_j(t)\psi
_j(s)$, and let $\mathbf{M}_k$ denote
the $S\times S$ matrix with elements $\sigma^{[k]}(\tau_r,\tau_s)$,
$r,s=1,\ldots,S$. Furthermore,
let $\lambda_{\mathrm{min}}(\mathbf{M}_k)$ denote the smallest eigenvalue of
the matrix $\mathbf{M}_k$.

%as4 #&#
\begin{assumption} \label{assum4}
(a) $\sup_t \sup_j \psi_j(t)^2\leq C_\psi$ for some $C_\psi<\infty$.\vspace*{-6pt}
\begin{longlist}[(b)]
\item[(b)] There exists some $0< C_1 <\infty$ such that $\lambda_j\leq
C_1 j^{-\mu}$ for all $j$.
\item[(c)] There
exists some $0<D<\infty$ such that $\lambda_{\mathrm{min}}(\mathbf
{M}_k)\geq D k^{-\mu+1}$ for all $k$.
\end{longlist}
\end{assumption}
Condition (a) is, for example, satisfied if $\psi_1,\psi_2,\ldots$
correspond to a Fourier-type basis.
Note that Assumption~\ref{assum3}(a) already implies that $\lambda_j$
must not be less
than a constant multiple of $j^{-\mu}$, and thus condition (b) requires
that $j^{-\mu}$ is also an upper bound for the rate of convergence of
$\lambda_j$. This in turn implies that $\sum_{j=k+1}^\infty\lambda
_j\leq C_2 k^{-\mu+1}$ as well as $|\sigma^{[k]}(t,s)| \leq C_2C_\psi^2
k^{-\mu+1}$
for some $ C_2<\infty$ and all $k$. Condition (c) therefore only
introduces an additional regularity condition on the matrix $\mathbf
{M}_k$. For the Brownian motion discussed in Section~\ref{sec3}, it is easily seen
that these requirements are necessarily fulfilled with $\mu=2$.

We now obtain the following theorem.

%th5 #&#
\begin{thmm} \label{parestGEN}
Under our setup and Assumptions \ref{assum2}--\ref{assum4} suppose that
$\widehat{S}=S$ and that
estimators $\widehat{\tau}_r$ satisfy (\ref{lem4eq1}) as well as (\ref
{lem4eq3}) for all $r=1,\ldots,S$.
If additionally
$k=O(n^{1/(\mu+2\nu)})$ and $n^{1/(\mu+2\nu)}=O(k)$ as $n\rightarrow
\infty$, then
%
%e5.4 #&#
%e5.5 #&#
\begin{eqnarray}
\Vert\widehat{\bolds{\beta}}-\bolds{\beta}\Vert _2^2&=&O_p
\bigl(n^{-2\nu/(\mu+2\nu)}\bigr), \label{thmpGENeq1}
\\
\int_a^b \bigl(\widehat{\beta}(t)-\beta(t)
\bigr)^2\,dt & =&O_p\bigl(n^{-(2\nu-1)/(\mu+2\nu)}\bigr).
\label{thmpGENeq2}
\end{eqnarray}
\end{thmm}

In the presence of points of impact the slope function $\beta(t)$ can
thus be estimated with the same
rate of convergence as in the standard model with $S=0$. The estimators
$\widehat{\beta}_r$ of $\beta_r$, $r=1,\ldots,S$, achieve a slightly
faster rate of convergence.

%t1 #&#
\begin{table}
\tabcolsep=0pt
\caption{Estimation errors for different sample sizes for the
simulation study
(OU-process, $\tau_1 =0.25$, $\tau_2=0.75$, $\beta_1= 2$, $\beta_2=
1$).
The column containing the estimate $\widehat{P}(\widehat{S}=S)$
contains two numbers: the estimate derived from the BIC followed by its
value derived from the cut-off procedure}\label{simtab}
\begin{tabular*}{\textwidth}{@{\extracolsep{\fill}}lccccccccccc@{}}
\hline
\multicolumn{2}{@{}l}{\textbf{Sample sizes}} & \multicolumn{10}{c@{}}{\textbf{Parameter
estimates}} \\[-6pt]
\multicolumn{2}{@{}l}{\hrulefill} & \multicolumn{10}{c@{}}{\hrulefill} \\
$\bolds{p}$ & $\bolds{n}$ & $\bolds{|\widehat{\tau}_1-\tau_1|}$ & $\bolds{|\widehat{\tau}_2-\tau_2|}$ &
$\bolds{|\widehat{\beta}_1-\beta_1|}$ & $\bolds{|\widehat{\beta}_2-\beta_2|}$ &$\bolds{\widehat
{S}} $ & $\bolds{\widehat{P}(\widehat{S}=S)}$ & $\bolds{\widehat{k} }$ & $\bolds{\smallint
(\widehat{\beta}-\beta)^2}$ %{$\langle\widehat{\beta}-\beta,\widehat{
%\beta}-\beta\rangle$}}
& \textbf{MSE} & \multicolumn{1}{c@{}}{$\bolds{\widehat{\kappa}}$} \\
\hline
\multicolumn{12}{c}{Simulation results if $\beta(t) \equiv0$} \\
\phantom{0.}1001 & \phantom{00}50 & 0.0130 & 0.0357 & 0.393 & 0.353 & 1.74 & 0.65/0.34 & 1.33
& 6.82 & 1.21 & 0.89 \\
& \phantom{0}100 & 0.0069 & 0.0226 & 0.274 & 0.249 & 1.96 & 0.77/0.40 & 1.05 &
3.43 & 1.21 & 0.94 \\
& \phantom{0}250 & 0.0027 & 0.0099 & 0.129 & 0.145 & 2.14 & 0.83/0.61 & 0.67 &
1.11 & 1.13 & 0.97 \\
& \phantom{0}500 & 0.0012 & 0.0061 & 0.070 & 0.097 & 2.15 & 0.86/0.73 & 0.45 &
0.51 & 1.08 & 0.98 \\
& 5000 & 0.0000 & 0.0004 & 0.012 & 0.012 & 2.04 & 0.96/0.98 & 0.03 &
0.00 & 1.00 & 1.00 \\
20,001& \phantom{00}50 & 0.0118 & 0.0333 & 0.393 & 0.350 & 1.71 & 0.64/0.35 & 1.78
& 6.91 & 1.19 & 0.89 \\
& \phantom{0}100 & 0.0068 & 0.0246 & 0.279 & 0.276 & 1.94 & 0.76/0.46 & 1.46 &
3.81 & 1.19 & 0.94 \\
& \phantom{0}250 & 0.0025 & 0.0108 & 0.121 & 0.144 & 2.15 & 0.83/0.62 & 0.74 &
1.02 & 1.12 & 0.97 \\
& \phantom{0}500 & 0.0013 & 0.0063 & 0.064 & 0.092 & 2.14 & 0.88/0.75 & 0.48 &
0.40 & 1.08 & 0.98 \\
& 5000 & 0.0001 & 0.0005 & 0.013 & 0.012 & 2.06 & 0.94/0.94 & 0.04 &
0.00 & 1.01 & 1.00 \\[3pt]
\multicolumn{12}{c}{Simulation results if $\beta(t) \neq0$} \\
\phantom{0.}1001 & \phantom{00}50 & 0.0150 & 0.0423 & 0.465 & 0.499 & 1.54 & 0.49/0.30 & 2.10
& 10.82 & 1.27 & 0.88 \\
& \phantom{0}100 & 0.0097 & 0.0317 & 0.376 & 0.400 & 1.86 & 0.63/0.34 & 2.06 &
5.93 & 1.27 & 0.94 \\
& \phantom{0}250 & 0.0039 & 0.0151 & 0.206 & 0.234 & 2.25 & 0.68/0.46 & 1.83 &
2.21 & 1.17 & 0.97 \\
& \phantom{0}500 & 0.0015 & 0.0083 & 0.107 & 0.164 & 2.30 & 0.72/0.59 & 1.69 &
0.90 & 1.10 & 0.99 \\
& 5000 & 0.0000 & 0.0006 & 0.036 & 0.027 & 2.25 & 0.79/0.97 & 2.01 &
0.05 & 1.01 & 1.00 \\
20,001& \phantom{00}50 & 0.0166 & 0.0399 & 0.467 & 0.465 & 1.52 & 0.47/0.29 & 2.14
& 11.19 & 1.29 & 0.89 \\
& \phantom{0}100 & 0.0099 & 0.0286 & 0.370 & 0.378 & 1.90 & 0.64/0.36 & 2.08 &
5.95 & 1.26 & 0.94 \\
& \phantom{0}250 & 0.0037 & 0.0171 & 0.185 & 0.263 & 2.27 & 0.67/0.49 & 1.90 &
2.19 & 1.15 & 0.97 \\
& \phantom{0}500 & 0.0018 & 0.0104 & 0.118 & 0.177 & 2.32 & 0.71/0.62 & 1.78 &
1.11 & 1.11 & 0.99 \\
& 5000 & 0.0002 & 0.0007 & 0.038 & 0.028 & 2.23 & 0.82/0.95 & 2.03 &
0.05 & 1.02 & 1.00 \\
\hline
\end{tabular*}
\end{table}

%%%%%%%%%%%%%%%%%%%%%%%%%%%%%
%s6 #&#
\section{Simulation study}\label{sec6}
%%%%%%%%%%%%%%%%%%%%%%%%%%%%%
We proceed by studying the finite sample performance of our estimation
procedure described in the preceding sections.
For different values of $n$, $p$, observations $(X_i,Y_i)$ are
generated according to the points of impact model \eqref{impact-model}
where $\varepsilon_i \sim N(0,1)$ are independent error terms. The
algorithms are implemented in R, and all tables are based on 1000
repetitions of the simulation experiments.
The corresponding R-code can be obtained from the authors upon request.

The data $X_1,\ldots,X_n$ are generated as independent
Ornstein--Uhlenbeck processes ($\kappa=1$) with parameters $\theta=5$
and $\sigma_u=3.5$ at $p$ equidistant grid points over the interval
$[0,1]$. Simulated trajectories are determined by using exact updating
formulas as proposed by \citet{Gillespie1996}.
The simulation study is based on $S=2$ points of impact located at
$\tau_1 = 0.25$ and $\tau_2 = 0.75$ with corresponding coefficients
$\beta_1 = 2$ as well as
$\beta_2 = 1$.
Results are reported in Table~\ref{simtab},
where the upper part of the table refers to the situation with $\beta
(t)\equiv0$,
while the lower part represents a model with $\beta(t) = 3.5t^3-5.5t^2+3t+0.5$.

In both cases, estimation of the points of impact relies on setting
$\delta=C \frac{1}{\sqrt{n}}$ for $C=1$, but similar results could be
obtained for a wide range of values $C$.
The results are then obtained by performing best subset selection with
the BIC-criterion via the R package bestglm on the augmented model
\eqref{augmentedmodel}
%
%e6.1 #&#
\begin{equation}
\label{augmentedmodel-2} Y_i\approx\sum_{r=1}^k
\alpha_r \langle X_i,\widehat{\psi}_r
\rangle +\sum_{r=1}^{\widetilde{S}}\beta_r
X_i(\tilde\tau_r)+\varepsilon_i^*.
\end{equation}
Here, $\widetilde{S}$ is the number of all possible candidates for the
points of impact and $k$ is initially set to $6$ principal components,
but tendencies remain unchanged for a broad range of values $k$.

For different sample sizes $n$ and $p$, Table~\ref{simtab} provides the
average absolute errors of our estimates, the frequency of $\widehat
{S}=S$, as well as average values of $\widehat{S}$, $\widehat{k}$, the
prediction error $\mathrm{MSE} = \frac{1}{n}\sum_{i=1}^n(\widehat
{y}_i-y_i)^2$ and $\widehat{\kappa}$. The column containing $\widehat
{P}(\widehat{S}=S)$ consists of two values. The first one being the
frequency of $\widehat{S}=S$ resulting from the BIC. For the second
one, $S$ was estimated by the cut-off
procedure using $\lambda=2\sqrt{\widehat{\operatorname{Var}}(Y)/n \log
(\frac{b-a}{\delta} )}$,
where $\widehat{\operatorname{Var}}(Y)$ denotes the estimated sample variance of $Y_i$.
The cut-off criterion yields very reliable estimates $\widehat{S}$ of
$S$ for $n=5000$, but showed a clear tendency to underestimate $S$ for
smaller sample sizes. The BIC-criterion however proves to possess a
much superior behavior in this regards for small $n$ but is
outperformed by the cut-off criterion for $n=5000$ in the case $\beta
(t) \neq0$.

In order to match $\{\widehat{\tau}_s\}_{s=1,\ldots,\widehat{S}}$ and $\{
\tau_r\}_{r=1,2}$ the interval $ [0,1 ]$ is divided into
$I_1= [0,\frac{1}{2}(\tau_{1}+\tau_{2}) ]$ and
$I_2 =  [\frac{1}{2}(\tau_{1}+\tau_{2}),1 ]$.
The estimate $\widehat{\tau}_s$ in interval $I_r$ with the minimal
distance to $\tau_r$ is then used as an estimate for $\tau_r$. No point
of impact candidate in interval $I_r$ results in an ``unmatched'' $\tau
_r$, $r=1,\ldots,S$ and a missing value when computing averages.

The table shows that estimates of points of impact are generally quite
accurate even for smaller sample sizes. The error decreases rapidly as
$n$ increases, and this improvement is essentially independent of $p$.
As expected, since $\beta_2<\beta_1$, the error of the absolute
distance between the second point of impact and its estimate is larger
than the error for the first point of impact.

Moreover, due to the common effect of the trajectory $X_i(\cdot)$ on
$Y_i$, the overall estimation
error in the case where $\beta(t)\neq0$ is slightly higher than in the
first case. At a first glance, one may be puzzled by the fact that for
$n=5000$ and $p=1001$ the average error $|\widehat{\tau}_r-\tau_r|$
is considerably smaller than the distance $\frac{1}{p-1}=\frac
{1}{1000}$ between two adjacent grid points. But note that our
simulation design implies that $\tau_r\in\{t_j| j=1,\ldots,p\}$,
$r=1,\ldots,S$, for $p=1001$ as well as $p=20\mbox{,}001$. For medium to large
sample sizes, there is thus a fairly high probability that $\widehat
{\tau}_r=\tau_r$. The case $p=1001$ particularly profits from this situation.
Finally, it can be seen that estimates for $\widehat{\kappa}$ tend to
slightly underestimate the true value $\kappa=1 $ for small values of $n$.

%%%%%%%%%%%%%%%%%%%%%%%%%%%%%
%s7 #&#
\section{Applications to real data}\label{sec7}
%%%%%%%%%%%%%%%%%%%%%%%%%%%%%
In this section, the algorithm from Section~\ref{SEC:4} is applied to a dataset
consisting of Canadian weather data. In this dataset, we relate the
mean relative humidity to hourly temperature data.
In the Supplementary Appendix A [Kneip, Poss and Sarda (\citeyear{KneipPossSarda2015})], a further
application can be found. We there analyze spectral data which play an
important role in spectrophotometry and different applied scientific fields.

In both examples, the algorithm is applied to centered observations and
the estimation procedure from Section~\ref{SEC:4} is modified by eliminating all
points in an interval of size $\delta|\log\delta|$ around a point of
impact candidate $\widehat{\tau}_j$, which is still sufficient to
establish assertion (\ref{thmm3eq1}).

After estimating $\widetilde{S}$ possible candidates for the points of
impact, the approximate model \eqref{augmentedmodel},
\[
Y_i\approx\sum_{r=1}^k
\alpha_r \langle X_i,\widehat{\psi}_r
\rangle +\sum_{r=1}^{\widetilde{S}}\beta_r
X_i(\tilde\tau_r)+\varepsilon_i^*,
\]
is used, where initially $k=6$ is chosen.
Over a fine grid of different values of~$\delta$, points of impact and
principal components are selected simultaneously by best subset
selection with the BIC-criterion and the model corresponding to the
minimal BIC is then chosen. The maximum number of variables selected by
the BIC-criterion is set to $6$ and all curves have been transformed to
be observed over $[0,1]$ when applying the algorithm from Section~\ref{SEC:4}.
The performance of the model is then measured by means of a
cross-validated prediction error.

In the Canadian weather dataset, the hourly mean temperature and
relative humidity from the $15$ closest weather stations in an area
around $100$ km from Montreal was obtained for each of the $31$ days in
December $2013$. The data was compiled from \url{http://climate.weather.gc.ca}.
Weather stations with more than ten missing observations on the
temperature or relative humidity were discarded from the dataset. The
remaining stations had their nonavailable observations replaced by the
mean of their closest observed predecessor and successor. After
preprocessing a total of $n=13$ weather stations remained and for each
station $p=744$ equidistant hourly observations of the temperature were
observed.
The response variable $Y_i$ was taken to be the mean over all observed
values of the relatively humidity at station $i$.

\begin{figure}

\includegraphics{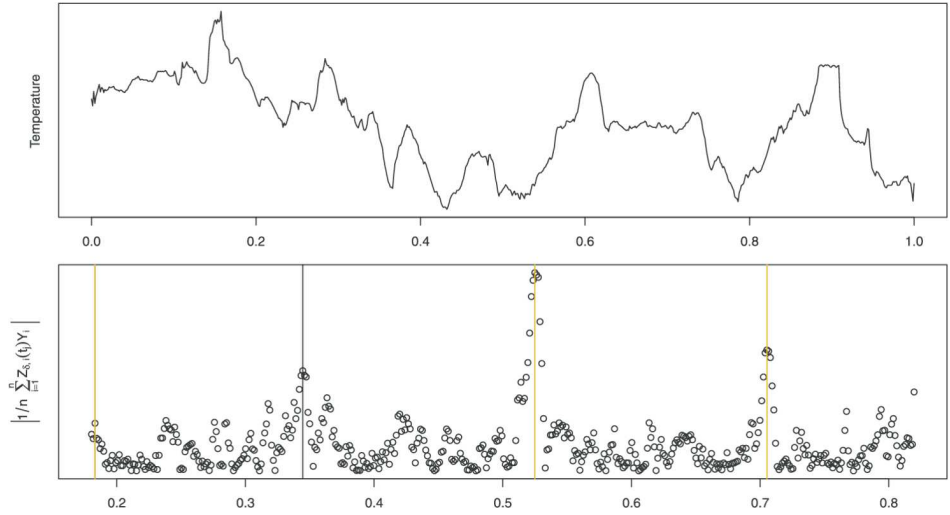}

\caption{The upper panel of this figure shows a
trajectory from the observed temperature curves of the Canadian weather
data. The lower panel shows {$|\frac{1}{n}\sum_{i=1}^n
Z_{\delta,i}(t_j)Y_i|$} during the selection procedure. Locations of
selected points of impact in the augmented model are indicated by
orange lines. The location of the remaining candidate is displayed by a
black line.}
\label{fig:Weather2}
\end{figure}

A cross-validated prediction error was calculated for three competing
regression models based on \eqref{augmentedmodel}.
In the first model, the mean relative humidity for each station was
explained by using the approximate model which combines the points of
impacts with a functional part.
The second and third model describe the cases $k = 0$ and $\widetilde
{S}=0$ in the approximate model, consisting only of points of impact
and the functional part, respectively.
For the first two models, points of impact were determined by
considering a total of $146$ equidistant values of $\delta$ between
$0.10$ and $0.49$. In all models BIC was used to approximate the
optimal values of the respective tuning parameters $\delta$, $S$ and/or
$k$ in a first step. The mean squared prediction error $\mathit{MSPE}=\frac
{1}{n}\sum_{i=1}^n(y_i - \hat{y}_i)^2$ was then calculated by means of
a leave one out cross-validation based on the chosen points of impact
and/or principal components from the first step. Additionally, the
median of $(y_i - \hat{y}_i)^2$, $i=1,\ldots,n$, has been calculated as
a more robust measure of the error.
 Depicted in the upper panel of Figure~\ref{fig:Weather2} is
the observed temperature trajectory for the weather station
``McTavish,'' showing a rather rough process.
The lower panel of this figure shows $|\frac{1}{n}\sum_{i=1}^n Z_{\delta,i}(t_j)Y_i|$ for the optimal value of $\delta=
0.18$ as obtained by the best model fit of the approximate model. While
orange lines represent the locations of the points of impact which were
actually selected with the help of the BIC-criterion, the location of
the remaining candidates are indicated by black vertical lines.

%f3 #&#
%
%t2 #&#
\begin{table}[b]
\caption{Estimated number of principal components $k$, points of impact
$S$,
prediction error and the median of $(y_i-\hat{y}_i)^2$ for the
Canadian weather data}
\label{app.tab:weather}
\begin{tabular*}{\textwidth}{@{\extracolsep{\fill}}lcccc@{}}
\hline
\textbf{Model} & $\bolds{\widehat{k}}$ & $\bolds{\widehat{S}}$ & $\mathbf{MSPE}$ & $\bolds{\operatorname{median}((y-\hat{y})^2)}$\\
\hline
Augmented & 3 & 3 & 2.314& 0.251\\% (12.86) \\
Points of impact & 0 & 3 & 1.714& 0.974\\ %(14.75)\\
FLR & 6 & 0 & 5.346& 1.269\\ %(11.83) \\
\hline
\end{tabular*}
\end{table}

Table~\ref{app.tab:weather} provides the empirical results when fitting
the three competing models. In terms of the prediction error, it can
clearly be seen from the table that the frequently applied functional
linear regression model is outperformed by the model consisting solely
of points of impact as well as the augmented (approximate) model. This
impression is supported by the last column of the table which gives the
median value of $(y_i - \hat{y}_i)^2$, showing additionally that,
typically, the augmented model performs even better than the plain
points of impact model.

An estimate $\widehat{\kappa}=0.14$ for $\kappa$ was obtained for
$\delta\approx0.3$, that is, the midpoint of the chosen values of
$\delta$. The estimated value of $\kappa=0.14$ corresponds to rather
rough sample paths as shown in the upper plot of Figure~\ref{fig:Weather2}.

In view of the small sample size results have to be interpreted with
care, and we therefore do not claim that this application provides
important substantial insights. Its main purpose is to serve as
illustration for classes of problems where our approach may be of
potential importance. It clearly shows that some relevant processes
observed in practice are nonsmooth.
With contemporary technical tools temperatures can be measured very
accurately, leading to a negligible measurement error. But
temperatures, especially in Canada, can vary rapidly over time. The
rough sample paths thus must be interpreted as an intrinsic feature of
temperature processes and cannot be explained by any type of ``error.''

\begin{appendix}
%%%%%%%%%%%%%%%%%%%%%%%%%%%%%
%s8 #&#
\section*{Appendix: Proofs of theorems}\label{app}

This appendix provides the
proofs of some of the main results. Remaining proofs can be
found in the supplementary material. Some of them rely on results from
\citeauthor{vandeGeerLederer2012} (\citeyear{vandeGeerLederer2012}),
\citeauthor{vanderVaartWellner2000} (\citeyear{vanderVaartWellner2000}) as well as
\citeauthor{ZhouLaffertyWassermann2008} (\citeyear{ZhouLaffertyWassermann2008}).

%%%%%%%%%%%%%%%%%%%%%%%%%%%%%
\begin{pf*}{Proof of Theorem~\ref{thmident}} Set $\beta_r:=0$ for
$r=S+1,\ldots,S^*$, and consider an arbitrary $j\in\{1,\ldots,S^*\}$.
Choose $0<\epsilon<\min_{r,s\in\{1,\ldots,S^*\},r\neq s} |\tau_r-\tau
_s|$ small enough such that conditions
(i)--(iv) of Definition~\ref{de1} are satisfied. Using (\ref{decompeps}), we
obtain a decomposition into two
uncorrelated components $X_{\epsilon,\tau_j}(\cdot)$ and $\zeta
_{\epsilon,\tau_j}(X) f_{\epsilon,\tau_j}(\cdot)$:
\begin{eqnarray*}
&&\mathbb{E} \Biggl( \Biggl(\int_a^b\bigl(
\beta(t)-\beta^*(t)\bigr)X(t)\,dt+\sum_{r=1}^{S^*}
\bigl(\beta_r-\beta_r^*\bigr)X(\tau_r)
\Biggr)^2 \Biggr)
\\
&&\qquad= \mathbb{E} \Biggl( \Biggl(\int_a^b\bigl(
\beta(t)-\beta^*(t)\bigr)X_{\epsilon,\tau_j}(t)\,dt+ \sum
_{r=1}^{S^*}\bigl(\beta_r-
\beta_r^*\bigr)X_{\epsilon,\tau_j}(\tau_r)
\Biggr)^2 \Biggr)
\\
&&\quad\qquad{} +\mathbb{E} \Biggl( \Biggl(\int_a^b\bigl(
\beta(t)-\beta^*(t)\bigr)\zeta _{\epsilon,\tau_j}(X) f_{\epsilon,\tau_j}(t)\,dt
\\
&&\qquad\quad{}+ \sum
_{r=1}^{S^*}\bigl(\beta_r-
\beta_r^*\bigr)\zeta_{\epsilon,\tau_j}(X) f_{\epsilon,\tau_j}(
\tau_r) \Biggr)^2 \Biggr)
\\
&&\qquad\geq \mathbb{E} \biggl( \biggl(\int_a^b
\bigl(\beta(t)-\beta^*(t)\bigr)\zeta _{\epsilon,\tau_j}(X)
 f_{\epsilon,\tau_j}(t)\,dt\\
 &&\qquad\quad{}+
\sum_{r\neq j}\bigl(\beta_r-
\beta_r^*\bigr)\zeta_{\epsilon,\tau_j}(X)f_{\epsilon
,\tau_j}(
\tau_r) + \bigl(\beta_j-\beta_j^*\bigr)
\zeta_{\epsilon,\tau_j}(X) f_{\epsilon,\tau_j}(\tau_j)
\biggr)^2 \biggr)
\\
&&\qquad\geq 2 \operatorname{var}\bigl(\zeta_{\epsilon,\tau_j}(X)\bigr) \bigl(
\beta_j-\beta_j^*\bigr) f_{\epsilon,\tau_j}(
\tau_j)\\
&&\qquad\quad{}\times \biggl( \int_a^b\bigl(
\beta(t)-\beta^*(t)\bigr) f_{\epsilon,\tau_j}(t)\,dt+\sum
_{r\neq
j}\bigl(\beta_r-\beta_r^*\bigr)
f_{\epsilon,\tau_j}(\tau_r) \biggr)
\\
&&\qquad\quad{} + \operatorname{var}\bigl(\zeta_{\epsilon,\tau_j}(X)\bigr) \bigl(
\beta_j-\beta _j^*\bigr)^2f_{\epsilon,\tau_j}(
\tau_j)^2.
\end{eqnarray*}

By condition (iv), we have
\[
\biggl|\sum_{r\neq j}\bigl(\beta_r-
\beta_r^*\bigr) f_{\epsilon,\tau_j}(\tau_r)\biggr| \leq \epsilon
S^* \max_{r\neq j}\bigl|\beta_r-\beta_r^*\bigr|
\bigl|f_{\epsilon,\tau
_j}(\tau_j)\bigr|,
\]
while boundedness of $\beta(\cdot)$ and $\beta^*(\cdot)$ implies that
there exits a constant
$0\leq D<\infty$ such that for all sufficiently small $\epsilon>0$
\begin{eqnarray*}
\biggl|\int_a^b\bigl(\beta(t)-\beta^*(t)\bigr)
f_{\epsilon,\tau_j}(t)\,dt\biggr| &\leq&\epsilon\int_{[a,b]\setminus[\tau_j-\epsilon,\tau_j+\epsilon]} D
\bigl|f_{\epsilon,\tau_j}(\tau_j)\bigr|\,dt \\
&&{}+\int_{\tau_j-\epsilon}^{\tau_j+\epsilon}
(1+\epsilon)D \bigl|f_{\epsilon
,\tau_j}(\tau_j)\bigr|\,dt
\\
&\leq&\epsilon\bigl(b-a+2(1+\epsilon)\bigr)D \bigl|f_{\epsilon,\tau_j}(
\tau_j)\bigr|.
\end{eqnarray*}
When combining these inequalities, we can conclude that for all
sufficiently small $\epsilon$ we have $\mathbb{E}(\int_a^b(\beta
(t)-\beta^*(t))X(t)\,dt+\sum_{r=1}^{S^*}(\beta_r-\beta_r^*)X(\tau_r))^2>0$
if $\beta_j-\beta_j^*\neq0$. Since $j\in\{1,\ldots,S^*\}$ is
arbitrary, the assertion of the theorem is an immediate consequence.
\end{pf*}

\begin{pf*}{Proof of Theorem~\ref{thmkarlov}} Choose some arbitrary $t\in
(a,b)$ and some $0<\epsilon<1$ with $\epsilon\leq\epsilon_t$. By
assumption, there exists a $k\in\mathbb{N}$ as well as some $f\in{\cal
C}(t,\epsilon,[a,b])$ such that
$|\langle f,\psi_r \rangle|>0$ for some $r\in\{1,\ldots,k\}$ and $\sup_{s\in[a,b]}|f_k(s)-f(s)|\leq\epsilon/3$, where $f_k(s)=\sum_{r=1}^k
\langle f,\psi_r \rangle\psi_r(s)$. The definition of ${\cal
C}(t,\epsilon,[a,b])$ then implies
that $f_k(t)\geq1-\epsilon/3$ as well as
%
%e8.1 #&#
\begin{eqnarray}
\label{karlovp1} \sup_{s\in[a,b]} \bigl|f_k(s)\bigr|&\leq&1+
\frac{\epsilon}{3}\leq(1+\epsilon) \biggl(1-\frac{\epsilon}{3}\biggr)\leq (1+
\epsilon)f_k(t),
\nonumber
\\[-8pt]
\\[-8pt]
\nonumber
\sup_{s\in[a,b],s\notin[t-\epsilon,t+\epsilon]} \bigl|f_k(s)\bigr|&\leq&\frac
{\epsilon}{3}\leq
\epsilon\biggl(1-\frac{\epsilon}{3}\biggr)\leq\epsilon f_k(t).
\end{eqnarray}
Now define the functional $\zeta_{\epsilon,t}$ by $\zeta_{\epsilon
,t}(X):=\sum_{r=1}^k \frac{\langle f,\psi_r \rangle}{\lambda_r} \langle
X,\psi_r \rangle$. Recall that
the coefficients $\langle X,\psi_r \rangle$ are uncorrelated and
$\operatorname{var}(\langle X,\psi_r \rangle)=\lambda_r$. By (\ref{karlov}), we
obtain
\begin{eqnarray*}
f_{\epsilon,t}(s)&:=&\frac{\mathbb{E}(X(s)\zeta_{\epsilon
,t}(X))}{\operatorname{var}(\zeta_{\epsilon,t}(X))} \\
&=&
\frac{\mathbb{E} ((\sum_{j=1}^\infty\langle X,\psi_j
\rangle\psi_j(s))(\sum_{r=1}^k ({\langle f,\psi_r \rangle}/{\lambda
_r}) \langle X,\psi_r\rangle) )}{\operatorname{var}(\zeta_{\epsilon,t}(X))}
\\
&=&\frac{\sum_{r=1}^k \langle f,\psi_r \rangle\psi_r(s)}{\operatorname{var}(\zeta
_{\epsilon,t}(X))} =\frac{f_k(s)}{\operatorname{var}(\zeta_{\epsilon,t}(X))}.
\end{eqnarray*}
Furthermore,
$\operatorname{var}(\zeta_{\epsilon,t}(X))=\sum_{r=1}^k \frac{\langle f,\psi_r \rangle
^2}{\lambda_r}>0$, and it thus follows from
(\ref{karlovp1}) that the functional $\zeta(t,X)$ satisfies conditions
(i)--(iv) of Definition~\ref{de1}.
Since $t\in(a,b)$ and $\epsilon$ are arbitrary, $X$ thus possesses
specific local variation.
\end{pf*}

\begin{pf*}{Proof of Theorem~\ref{thmcharX}} First note that Assumption~\ref
{assum1} implies that the absolute values of all
first and second-order partial derivatives of $\omega(t,s,z)$ are
uniformly bounded\vspace*{1pt}
by some constant $M<\infty$ for all $(t,s,z)$ in the compact subset
$[a,b]^2\times[0,b-a]$
of $\Omega$.

By definition of $Z_\delta$, it thus follows from a Taylor expansion
of $\omega$ that
for $t\in(a,b)$, any sufficiently small $\delta>0$ and some constant
$M_1<\infty$
%
%e8.2 #&#
\begin{eqnarray}
\label{Zdelta1} \mathbb{E}\bigl(X(t)Z_{\delta}(X,t)\bigr)&=&\sigma(t,t)-
\frac{1}{2}\sigma (t,t-\delta) -\frac{1}{2}\sigma(t,t+\delta)
\nonumber
\\
&=& \omega(t,t,0)-\frac{1}{2}\omega\bigl(t,t-\delta,\delta^\kappa
\bigr) -\frac{1}{2}\omega\bigl(t,t+\delta,\delta^\kappa\bigr)
\\
&=& \delta^\kappa c(t) +R_{1;\delta,t} \qquad
\mbox{with } \sup
_{t\in[a+\delta,b-\delta]} |R_{1;\delta,t}|\leq M_1
\delta^{\min\{2\kappa,2\}}.\nonumber
\end{eqnarray}
For the variance of $Z_{\delta}(X,t)$, we obtain by similar arguments
%
%e8.3 #&#
\begin{eqnarray}
\label{Zdelta2} \operatorname{var}\bigl(Z_{\delta}(X,t)\bigr)&=& 2
\omega(t,t,0)-\omega\bigl(t,t-\delta,\delta^\kappa\bigr)-\omega\bigl(t,t+
\delta,\delta ^\kappa\bigr)\nonumber\\
&&{} - \frac{1}{2} \bigl(\omega(t,t,0)-\omega
\bigl(t+\delta,t-\delta,(2\delta )^\kappa\bigr) \bigr)
\nonumber
\\
&&{} - \frac{1}{4} \bigl(2\omega(t,t,0)-\omega(t-\delta,t-\delta,0) -
\omega(t+\delta,t+\delta,0) \bigr)
\\
&=&\delta^\kappa \biggl(2c(t)-\frac{2^\kappa}{2}c(t) \biggr)
+R_{2;\delta,t}\nonumber\\
\eqntext{ \mbox{with } \sup_{t\in[a+\delta,b-\delta]}
|R_{2;\delta
,t}|<M_2\delta^{\min\{2\kappa,2\}}}
\end{eqnarray}
for some constant $M_2<\infty$.
Moreover, for any $0<c<\infty$ Taylor expansions of $\omega$ yield that
for any sufficiently small $\delta>0$ and all $u\in[-c,c]$
%
%e8.4 #&#
%e8.5 #&#
\begin{eqnarray}
&&\mathbb{E} \bigl(X(t+u\delta)Z_{\delta}(X,t)\bigr)\nonumber\\
&&\qquad=\sigma(t+u\delta,t)-
\tfrac
{1}{2}\sigma(t+u\delta,t-\delta) -\tfrac{1}{2}\sigma(t+u\delta,t+
\delta)
\nonumber
\\
&&\qquad= \omega(t,t,0) -\tfrac{1}{2}\omega\bigl(t,t-\delta,\delta^\kappa
\bigr) -\tfrac{1}{2}\omega\bigl(t,t+\delta,\delta^\kappa\bigr)
\nonumber
\\
&&\qquad\quad{} -c(t)\delta^\kappa \bigl(|u|^\kappa -\tfrac{1}{2}
\bigl(|u+1|^\kappa-1\bigr)-\tfrac{1}{2}\bigl(|u-1|^\kappa-1
\bigr) \bigr) + R_{3;c,u,\delta,t} \label{Zdelta3}
\\
&&\qquad= -c(t)\delta^\kappa \bigl(|u|^\kappa -\tfrac{1}{2}|u+1|^\kappa-
\tfrac{1}{2}|u-1|^\kappa \bigr) + R_{4;c,u,\delta,t}, \label{Zdelta4}
\end{eqnarray}
where for some constants $M_{3,c}<\infty$ and $M_{4,c}<\infty$
\begin{eqnarray*}
\sup_{t\in[a+\delta,b-\delta]} R_{3;c,u,\delta,t}&\leq& M_{3,c}
\bigl(|u|^{1/2}\delta\bigr)^{\min\{2\kappa,2\}},\\
 \sup_{t\in[a+\delta,b-\delta]}
R_{4;c,u,\delta,t}&\leq& M_{4,c}\delta ^{\min\{2\kappa,2\}}
\end{eqnarray*}
hold for all $u\in[-c,c]$.
Finally, Assumption~\ref{assum1} implies that there exists a constant
$M_5<\infty$ such that
for all $s\in[a,b]$ with $|t-s|\geq\delta$\vspace*{-2pt}
%
%
%e8.6 #&#
\begin{eqnarray}
\label{Zdelta5}&& \bigl|\mathbb{E} \bigl(X(s)Z_{\delta}(X,t)\bigr)\bigr|\nonumber\\[-1pt]
&&\qquad= \bigl| \omega
\bigl(s,t,|s-t|^\kappa\bigr) -\tfrac{1}{2}\omega\bigl(s,t-
\delta,|s-t+\delta|^\kappa\bigr)\nonumber
\\[-8pt]
\\[-8pt]
&&\qquad\quad{} -\tfrac{1}{2}\omega\bigl(s,t+
\delta,|s-t-\delta|^\kappa\bigr)\bigr|\nonumber
\\[-1pt]
&&\qquad\leq %
\cases{ \displaystyle M_5 \frac{\delta^2}{|t-s|^{2-\kappa}}, &\quad $ \mbox{if }
\kappa\neq1,$\vspace*{2pt}
\cr
M_5 \delta^2, &\quad $\mbox{if }
\kappa= 1$.}\nonumber %
\end{eqnarray}

It follows from (\ref{Zdelta1}), (\ref{Zdelta4}) and (\ref{Zdelta5})
that for arbitrary $t\in(a,b)$ and any $\epsilon>0$ there exist a
$\delta_\epsilon>0$ as well as a constant
$a_\epsilon\geq1$ such that for all $\delta\leq\delta_\epsilon$\vspace*{-2pt}
\begin{eqnarray*}
\bigl|\mathbb{E} \bigl(X(s)Z_{\delta}(X,t)\bigr)\bigr|&\leq&(1+\epsilon)\mathbb{E}
\bigl(X(t)Z_{\delta}(X,t)\bigr)\qquad \mbox{for all } s\in[a,b], s\neq t,
\\
\bigl|\mathbb{E} \bigl(X(s)Z_{\delta}(X,t)\bigr)\bigr|& \leq&\epsilon\cdot\mathbb{E}
\bigl(X(t)Z_{\delta}(X,t)\bigr) \qquad\mbox{for all } s\in[a,b], |s-t|\geq
a_\epsilon \delta.
\end{eqnarray*}
Together with (\ref{Zdelta2}), the assertion of the theorem is an
immediate consequence.
\end{pf*}

\begin{pf*}{Proof of Theorem~\ref{parestGEN}}
Let $\hat\theta_{ij}:=\langle X_i,\widehat{\psi}_j\rangle$, $ \theta
_{ij}:=\langle X_i,\psi_j\rangle$, and
$\tilde\alpha_j:= \langle\beta,\widehat{\psi}_j\rangle$ for all $i,j$.
Using empirical eigenfunctions, we obtain $X_i=\sum_{j=1}^n \hat\theta
_{ij}\widehat{\psi}_j$ and $\int_a^b \beta(t)X_i(t)\,dt =\sum_{j=1}^n
\tilde\alpha_j\hat\theta_{ij}$. Therefore,\vspace*{-2pt}
%
%e8.7 #&#
\begin{equation}
Y_i=\sum_{j=1}^n \Biggl(
\tilde\alpha_j +\sum_{r=1}^S
\beta_r \widehat {\psi}_j(\tau_r) \Biggr)
\hat\theta_{ij}+\varepsilon_i, \label{parGEN-eq1}
\end{equation}
and for all possible values $b_1,\ldots,b_S$ and all $a_1,\ldots,a_k$\vspace*{-2pt}
%
%e8.8 #&#
\begin{eqnarray}\label{parGEN-eq2}
&&\sum_{j=1}^k a_j \hat
\theta_{ij}+\sum_{r=1}^S
b_r X_i(\widehat{\tau}_r)
\nonumber
\\[-9pt]
\\[-9pt]
\nonumber
&&\qquad = \sum
_{j=1}^k \Biggl( a_j +\sum
_{r=1}^S b_r \widehat{
\psi}_j(\widehat {\tau}_r) \Biggr) \hat
\theta_{ij}+\sum_{j=k+1}^n \sum
_{r=1}^Sb_r \widehat{
\psi}_j(\widehat{\tau}_r)\hat\theta_{ij}
\end{eqnarray}
for all $i=1,\ldots,n$. By definition, $\hat\lambda_j=\frac{1}{n}\sum_{i=1}^n \hat\theta_{ij}^2$, $j=1,\ldots,n$, and for $j\ne l$ the
coefficients $\hat\theta_{ij}$
and $\hat\theta_{il}$ are empirically uncorrelated, that is, $\sum_{i=1}^n \hat\theta_{ij}\hat\theta_{il}=0$. It follows that for any
given values
$b_1,\ldots,b_S$ the values $\hat\alpha(\mathbf{b})_j$, $j=1,\ldots
,k$, minimizing
$\sum_{i=1}^n  (
Y_i-\sum_{j=1}^k a_j \hat\theta_{ij} - \sum_{r=1}^{S}b_r X_i(\widehat
\tau_r) )^2$ over all $a_1,\ldots,a_k$ are given by
%
%e8.9 #&#
\begin{eqnarray}\label{parGEN-eq3}
\hat\alpha(\mathbf{b})_j=\tilde\alpha_j+\hat
\lambda_j^{-1} \frac{1}{n}\sum
_{i=1}^n\hat\theta_{ij}
\varepsilon_i+ \sum_{r=1}^S
\bigl(\beta_r \widehat{\psi}_j(\tau_r)-b_r
\widehat{\psi }_j(\widehat{\tau}_r)\bigr),
\nonumber
\\[-10pt]
\\[-10pt]
\eqntext{j=1,
\ldots,k. }
\end{eqnarray}
Note that $\tilde\alpha_j+\hat\lambda_j^{-1}
\frac{1}{n}\sum_{i=1}^n\hat\theta_{ij}\varepsilon_i$ is identical to
the estimate
of $\alpha_j$ to be obtained in a standard functional\vadjust{\goodbreak} linear regression
model with no points of impact.
Theorem~1 of \citet{HallHorowitz2007} thus implies that
%
%e8.10 #&#
\begin{eqnarray}\label{parGEN-eq4}
&&\int_a^b \Biggl(\beta(t)-\sum
_{j=1}^k\Biggl( \tilde\alpha_j+\hat
\lambda_j^{-1} \frac{1}{n}\sum
_{i=1}^n\hat\theta_{ij}
\varepsilon_i\Biggr)\widehat{\psi }_j(t)\,dt
\Biggr)^2\,dt
\nonumber
\\[-9pt]
\\[-9pt]
\nonumber
&&\qquad  =O_p\bigl(n^{-(2\nu-1)/(\mu+2\nu)}\bigr).
\end{eqnarray}
Further analysis requires to analyze the differences between $\theta
_{ij}, \psi_j$ and
their empirical counterparts
$\hat\theta_{ij}, \widehat{\psi}_j$.
By Assumptions \ref{assum2}--\ref{assum4} and
$k=O(n^{1/(\mu+2\nu)})$,
Theorems 1 and 2 together with equation (2.8) of \citet{HallHosseini2006} imply that
for any $q=1,2,3,\ldots$ there exists some $A_q,B_q<\infty$ such that
%
%e8.11 #&#
\begin{eqnarray}\label{parGEN-eq5}
E \bigl(|\lambda_j-\hat\lambda_j|^q
\bigr)&\leq& A_qn^{-q/2},
\nonumber
\\[-9pt]
\\[-9pt]
\nonumber
 \sup_t E
\bigl(\bigl|\widehat{\psi}_j(t)-\psi_j(t)\bigr| \bigr)&\leq&
B_qn^{-q/2}j^{q(\mu
+1)},\qquad  j=1,\ldots,k+1 \vspace*{-2pt}
\end{eqnarray}
for all sufficiently large $n$. Let $X_i^{[k]}:= X_i-\sum_{j=1}^k \hat
\theta_{ij}\widehat{\psi}_j$.
Recall that\break $\lambda_j=O(j^{-\mu})$ and note that by Assumptions \ref
{assum3} and \ref{assum4}, $n^{-1/2}n^{2/(\mu+2\nu)}=\break o(n^{(-\mu+1)/(\mu
+2\nu)})$, while $ n^{(-\mu+1)/(\mu+2\nu)}=O(\sigma^{[k]}(\tau_r,\tau_r))$.
By (\ref{parGEN-eq5}), we thus obtain for all $t,s\in[a,b]$
%
%e8.12 #&#
\begin{eqnarray}
\label{parGEN-eq6} &&\frac{1}{n} \sum_{i=1}^n
X_i^{[k]}(t)X_i^{[k]}(s)\nonumber\\[-1pt]
&&\qquad=
\frac{1}{n} \sum_{i=1}^n
X_i(t)X_i(s) -\sum_{j=1}^k
\hat\lambda_j \widehat{\psi}_j(t)\widehat{
\psi}_j(s)
\nonumber
\\[-1pt]
&&\qquad=\sigma(t,s)-\sum_{j=1}^k
\lambda_j \psi_j(t)\psi_j(s) +\sum
_{j=1}^k \lambda_j\bigl(
\psi_j(t)\psi_j(s)-\widehat {\psi}_j(t)
\widehat{\psi}_j(s)\bigr)
\nonumber
\\[-9pt]
\\[-9pt]
\nonumber
&&\qquad\quad{} +\sum_{j=1}^k (\lambda_j-
\hat\lambda_j)\widehat{\psi }_j(t)\widehat{
\psi}_j(s) +O_P\bigl(n^{-1/2}\bigr)
\nonumber
\\[-1pt]
&&\qquad=\sigma^{[k]}(t,s) +O_P\bigl(n^{-1/2}n^{2/(\mu+2\nu)}
\bigr)\nonumber\\[-1pt]
&&\qquad=\sigma^{[k]}(t,s) +o_P\bigl(n^{(-\mu+1)/(\mu+2\nu)}\bigr).\nonumber
\end{eqnarray}
At the same time, (\ref{lem4eq1}) leads to
%
%e8.13 #&#
\begin{eqnarray}
\label{parGEN-eq7}&& \frac{1}{n} \sum_{i=1}^n
\bigl(X_i^{[k]}(\tau_r)-X_i^{[k]}(
\widehat{\tau }_r)\bigr)^2\nonumber\\
&&\qquad= \frac{1}{n} \sum
_{i=1}^n \bigl(X_i(
\tau_r)-X_i(\widehat{\tau }_r)
\bigr)^2 -\sum_{j=1}^k \hat
\lambda_j \bigl(\widehat{\psi}_j(\tau_r)-
\widehat{\psi}_j(\widehat{\tau}_r)\bigr)^2
\\
&&\qquad\leq \frac{1}{n} \sum_{i=1}^n
\bigl(X_i(\tau_r)-X_i(\widehat{\tau
}_r)\bigr)^2=O_P\bigl(n^{-1}
\bigr)\nonumber
\end{eqnarray}
for all $r=1,\ldots,S$. Expressions (\ref{parGEN-eq6}) and (\ref
{parGEN-eq7}) together imply that for
all $r,s$
%
%e8.14 #&#
\begin{equation}
\frac{1}{n} \sum_{i=1}^n
X_i^{[k]}(\widehat{\tau}_r)
X_i^{[k]}(\widehat {\tau}_s)=
\sigma^{[k]}(\tau_r,\tau_s)+o_P
\bigl(n^{(-\mu+1)/(\mu+2\nu)}\bigr). \label{parGEN-eq8}
\end{equation}
Let $\mathbf{X}_i^{[k]}:=(X_i^{[k]}(\widehat{\tau}_1),\ldots
,X_i^{[k]}(\widehat{\tau}_S))^T$ and note
that by (\ref{parGEN-eq8}) we have\break  $\frac{1}{n} \sum_{i=1}^n \mathbf
{X}_i^{[k]}(\mathbf{X}_i^{[k]})^T$ $=\mathbf{M}_k+
o_P(n^{(-\mu+1)/(\mu+2\nu)})$.
By Assumption~\ref{assum4}(b), we can conclude that with probability
tending to 1 as $n\rightarrow\infty$ the matrix
$\frac{1}{n}\sum_{i=1}^n \mathbf{X}_i^{[k]}(\mathbf{X}_i^{[k]})^T$ is
invertible,
%
%e8.15 #&#
\begin{eqnarray} \label{parGEN-eq9}
&&n^{(-\mu+1)/(\mu+2\nu)}\Biggl(\frac{1}{n}\sum_{i=1}^n
\mathbf {X}_i^{[k]}\bigl(\mathbf{X}_i^{[k]}
\bigr)^T\Biggr)^{-1}
\nonumber
\\[-8pt]
\\[-8pt]
\nonumber
&&\qquad=n^{(-\mu+1)/(\mu+2\nu
)}(\mathbf{M}_k)^{-1}
+o_P(1)
\end{eqnarray}
and hence by
(\ref{parGEN-eq1})--(\ref{parGEN-eq3}) the least squares
estimator $\widehat{\bolds{\beta}}$ of $\bolds{\beta}$ can be
written in the form
%
%e8.16 #&#
\begin{eqnarray}\label{parGEN-eq10}
\widehat{\bolds{\beta}}&=&\Biggl(\frac{1}{n}\sum
_{i=1}^n \mathbf {X}_i^{[k]}
\bigl(\mathbf{X}_i^{[k]}\bigr)^T
\Biggr)^{-1}
\nonumber
\\[-8pt]
\\[-8pt]
\nonumber
&&{}\times\frac{1}{n}\sum_{i=1}^n
\mathbf{X}_i^{[k]} \Biggl(\sum
_{r=1}^S \beta_r X_i^{[k]}(
\tau_r) + \sum_{j=k+1}^n \tilde
\alpha_j \hat\theta_{ij}+\varepsilon_i
\Biggr).
\end{eqnarray}
By (\ref{parGEN-eq7}) and (\ref{parGEN-eq8}), we obtain
%
%e8.17 #&#
\begin{eqnarray} \label{parGEN-eq10a}
&&\frac{1}{n}\sum_{i=1}^n
\mathbf{X}_i^{[k]}\sum_{r=1}^S
\beta_r X_i^{[k]}(\tau_r)
\nonumber
\\[-8pt]
\\[-8pt]
\nonumber
&&\qquad=
\frac{1}{n}\sum_{i=1}^n
\mathbf{X}_i^{[k]}\bigl(\mathbf {X}_i^{[k]}
\bigr)^T\bolds{\beta}+ O_P\bigl(n^{(-\mu+1)/2(\mu+2\nu)}\cdot
n^{-1/2}\bigr).
\end{eqnarray}
The results of \citet{HallHorowitz2007} imply that $\sum_{j=k+1}^n
\tilde\alpha_j^2= \break O_P(n^{-(2\nu-1)/(\mu+2\nu)})$. The Cauchy--Schwarz
inequality thus leads to
%
%e8.18 #&#
\begin{eqnarray}
\label{parGEN-eq11}&& \Biggl|\frac{1}{n}\sum_{i=1}^n
X_i^{[k]}(\widehat{\tau}_r) \Biggl(\sum
_{j=k+1}^n \tilde\alpha_j \hat
\theta_{ij}\Biggr)\Biggr|\nonumber\\
&&\qquad =\Biggl|\sum_{j=k+1}^n
\tilde\alpha_j \hat\lambda_j \widehat{
\psi}_j(\widehat {\tau}_r)\Biggr|
\nonumber
\\[-8pt]
\\[-8pt]
\nonumber
&&\qquad\leq\sqrt{\sum_{j=k+1}^n \hat
\lambda_j \tilde\alpha_j^2}\sqrt{
\sum_{j=k+1}^n \hat\lambda_j
\widehat{\psi}_j(\widehat{\tau}_r)^2} \leq
\sqrt{\hat\lambda_{k+1} \sum_{j=k+1}^n
\tilde\alpha_j^2}\sqrt {\frac{1}{n}\sum
_{i=1}^n X_i^{[k]}(
\widehat{\tau}_r)^2}
\\
&&\qquad =O_P\bigl( n^{-(\mu+2\nu-1)/2(\mu+2\nu)}\cdot
n^{(-\mu+1)/2(\mu+2\nu)}\bigr)\nonumber
\end{eqnarray}
for all $r=1,\ldots,S$. Furthermore, $\widehat{\psi}_j(t)=
\hat\lambda_j^{-1}\frac{1}{n}\sum_{i=1}^n\hat\theta_{ij}X_i(t)$, and hence
the Cauchy--Schwarz inequality yields
%
%e8.19 #&#
\begin{eqnarray}\label{parGEN-eq11a}
\bigl|\widehat{\psi}_j(\tau_r)-\widehat{
\psi}_j(\widehat{\tau}_r)\bigr|&=& \Biggl|\hat\lambda_j^{-1}
\frac{1}{n}\sum_{i=1}^n \hat
\theta_{ij}\bigl(X_i(\tau_r)-X_i(
\widehat{\tau }_r)\bigr)\Biggr|
\nonumber
\\[-8pt]
\\[-8pt]
\nonumber
&\leq &\hat\lambda_j^{-1/2}
\sqrt{\frac{1}{n}\sum_{l=1}^n
\bigl(X_l(\tau _r)-X_l(\widehat{
\tau}_r)\bigr)^2}.
\end{eqnarray}
Now note that by the
independence of $\hat\theta_{ij}$ and $\varepsilon_i$ we have $\hat
\lambda_j^{-1/2}\frac{1}{n}\sum_{i=1}^n\hat\theta_{ij}\varepsilon
_i=O_P(n^{-1/2})$. By
(\ref{lem4eq3}), it therefore follows from (\ref{parGEN-eq11a})
that
\begin{eqnarray*}
&&\frac{1}{n}\sum_{i=1}^n
\bigl(X_i^{[k]}(\widehat{\tau}_r)-X_i^{[k]}(
\tau _r)\bigr)\varepsilon_i\\
&&\qquad=\frac{1}{n}\sum
_{i=1}^n \bigl(X_i(\widehat{
\tau}_r)-X_i(\tau_r)\bigr)
\varepsilon_i- \sum_{j=1}^k
\frac{1}{n}\sum_{i=1}^n\hat
\theta_{ij}\varepsilon_i \bigl(\widehat{
\psi}_j(\widehat{\tau}_r)-\widehat{\psi}_j(
\tau_r)\bigr)
\\
&&\qquad=O_P\bigl((k+1) n^{-1}\bigr)=O_P\bigl(
n^{-(\mu+2\nu-1)/(\mu+2\nu)}\bigr).
\end{eqnarray*}
Using (\ref{parGEN-eq6}), it is immediately seen that
$\frac{1}{n}\sum_{i=1}^n X_i^{[k]}(\tau_r)\varepsilon_i=\break O_P(n^{-1/2}
n^{(-\mu+1)/2(\mu+2\nu)})$. Consequently,
%
%e8.20 #&#
\begin{eqnarray} \label{parGEN-eq12}
\frac{1}{n}\sum_{i=1}^n
X_i^{[k]}(\widehat{\tau}_r)
\varepsilon_i&=& \frac{1}{n}\sum_{i=1}^n
X_i^{[k]}(\tau_r)\varepsilon_i +
\frac{1}{n}\sum_{i=1}^n
\bigl(X_i^{[k]}(\widehat{\tau}_r)-X_i^{[k]}(
\tau _r)\bigr)\varepsilon_i
\nonumber
\\[-8pt]
\\[-8pt]
\nonumber
&=&O_P
\bigl(n^{-1/2} n^{(-\mu+1)/2(\mu+2\nu)}\bigr).
\end{eqnarray}
By Assumption~\ref{assum4}(c), we can infer from (\ref{parGEN-eq9})
that the maximal
eigenvalue of the matrix $(\frac{1}{n}\sum_{i=1}^n \mathbf
{X}_i^{[k]}(\mathbf{X}_i^{[k]})^T)^{-1}$
can be bounded by $\lambda_{\mathrm{max}} ((\frac{1}{n}\sum_{i=1}^n \mathbf
{X}_i^{[k]}\times\break (\mathbf{X}_i^{[k]})^T)^{-1})=O_P(n^{(\mu-1)/(\mu+2\nu)})$.
It therefore follows
from (\ref{parGEN-eq10})--(\ref{parGEN-eq12}) that
\begin{eqnarray*}
\widehat{\bolds{\beta}}&=& \bolds{\beta}+ O_P\bigl(n^{(\mu-1)/(\mu+2\nu)}
\cdot n^{(-\mu+1)/2(\mu+2\nu)}\cdot n^{-(\mu+2\nu-1)/2(\mu+2\nu)}\bigr)\\
&= &\bolds{\beta}+O_P
\bigl( n^{-\nu/(\mu+2\nu)}\bigr).
\end{eqnarray*}
This proves (\ref{thmpGENeq1}). Using (\ref{parGEN-eq3}), it follows
that the
least squares estimators $\widehat{\alpha}_j$ of $\tilde\alpha_j$ are
given by
%
%e8.21 #&#
\begin{eqnarray}\label{parGEN-eq13}
\widehat{\alpha}_j&=&\tilde\alpha_j+\hat
\lambda_j^{-1} \frac{1}{n}\sum
_{i=1}^n\hat\theta_{ij}
\varepsilon_i+ \sum_{r=1}^S (
\beta_r-\widehat{\beta}_r) \widehat{
\psi}_j(\tau_r)
\nonumber
\\[-8pt]
\\[-8pt]
\nonumber
&&{}-\sum_{r=1}^S
\widehat{\beta}_r\bigl( \widehat{\psi}_j(\widehat{
\tau}_r)-\widehat{\psi}_j(\tau_r)\bigr),\qquad
j=1,\ldots,k .
\end{eqnarray}
But (\ref{parGEN-eq5}) and (\ref{thmpGENeq1}) imply that
%
%e8.22 #&#
\begin{eqnarray} \qquad\hspace*{6pt}\label{parGEN-eq14}
\sum_{j=1}^k \Biggl(\sum
_{r=1}^S (\beta_r-\widehat{
\beta}_r) \widehat{\psi }_j(\tau_r)
\Biggr)^2= O_P\bigl(k n^{-2\nu/(\mu+2\nu)}
\bigr)
=O_P\bigl(n^{-(2\nu-1)/(\mu+2\nu)}\bigr),\hspace*{-8pt}
\end{eqnarray}
while by (\ref{parGEN-eq5}) and (\ref{parGEN-eq11a})
\begin{eqnarray*}
\sum_{j=1}^k\bigl(\widehat{
\psi}_j(\tau_r)-\widehat{\psi}_j(\widehat{
\tau}_r)\bigr)^2 \leq\frac{k}{\lambda_k}
\frac{1}{n}\sum_{i=1}^n
\bigl(X_i(\tau_r)-X_i(\widehat{
\tau}_r)\bigr)^2 =O_P\bigl(n^{-(2\nu-1)/(\mu+2\nu)}
\bigr),
\end{eqnarray*}
and therefore
%
%e8.23 #&#
\begin{equation}
\sum_{j=1}^k \Biggl(\sum
_{r=1}^S\widehat{\beta}_r\bigl(
\widehat{\psi}_j(\widehat{\tau}_r)-\widehat{
\psi}_j(\tau_r)\bigr)\Biggr)^2=
O_P\bigl(n^{-(2\nu-1)/(\mu+2\nu)}\bigr) \label{parGEN-eq15}.
\end{equation}
Assertion (\ref{thmpGENeq2}) now is an immediate consequence of (\ref
{parGEN-eq4}) and
(\ref{parGEN-eq13})--(\ref{parGEN-eq15}).
\end{pf*}
\end{appendix}

\begin{supplement}%[id=suppA]
\stitle{Supplement to ``Functional linear regression with points of
impact''}
\slink[doi]{10.1214/15-AOS1323SUPP} %[doi,text={...}] - jei reikia suskaldyti doi
\sdatatype{.pdf}
\sfilename{aos1323\_supp.pdf}
\sdescription{The supplementary document by Kneip, Poss and Sarda (\citeyear{KneipPossSarda2015})
contains three Appendices. An application to NIR data can be found in
Appendix~A. In Appendix~B, it is shown that the eigenfunctions of a
Brownian motion satisfy assertion \ref{karlovp0} in Theorem~\ref
{thmkarlov}. Appendix~C provides the proofs of Theorem~\ref{poicon} and
Propositions~\ref{poiconkappa} and \ref{lem4}.}
\end{supplement}

% imsref loaded by akundreckaite, 2015-08-27 13:53:07

%
%\begin{appendix}
%\section{}
%\end{appendix}

%\section*{Acknowledgments}

%\begin{thebibliography}{99}
%\bibitem[\protect\citeauthoryear{}{}]{r1}
%\bibitem{r1}
%\end{thebibliography}

\printaddresses
\end{document}